\documentclass[a4paper,leqno]{article}

\usepackage[centertags]{amsmath}
\usepackage{amsfonts}
\usepackage{amssymb}
\usepackage{amsthm}

\def \Real{\mbox{\sl I\kern-.166em R}}            
    \def \Nat{\mbox{\sl I\kern-.166em N}}             
\def \Pr{\mbox{\sl I\kern-.166em P}}

    \def\refart#1#2#3#4{#1, #2, {\em #3\/}\ #4.}
    \def\book#1#2#3{#1, {\em #2,} #3.}

\newcommand{\red}{\protect\large\bf}

\newtheorem{Th}{Theorem}
\newtheorem{Cor}{Corollary}
\newtheorem{Pro}{Proposition}

\theoremstyle{definition}

\newtheorem{Conv}{CONVENTION}

\newtheorem{Nots}{Notations}
\newtheorem{Rk}{Remark}

\numberwithin{equation}{section}

\title{\bf\large 	HARMONIC SECTIONS OF TANGENT BUNDLES EQUIPPED WITH RIEMANNIAN $g$-NATURAL METRICS}

\author{M.T.K. ABBASSI\thanks{ This work was started during the
visit of the first author at the University of Lecce. The first
author would like to express his sincere thanks to the second
and third authors for their kind invitation and their
hospitality.
\newline
\indent ${}^{**}$Authors supported by funds of the University of
Lecce and M.I.U.R. (PRIN 2005).
} ,        
 $ $ G. CALVARUSO \footnotemark[7]  $ $ $ $ and D. PERRONE \footnotemark[7]
 }

\date{}
\begin{document}
\maketitle

\begin{abstract}
Let $(M,g)$ be a Riemannian manifold. When $M$ is compact and the tangent bundle $TM$ is equipped with the Sasaki metric $g^s$, the only vector fields which define harmonic maps from $(M,g)$ to $(TM,g^s)$, are the parallel ones. The Sasaki metric, and other well known Riemannian metrics on $TM$, are particular examples of $g$-natural metrics. We equip $TM$ with an arbitrary Riemannian $g$-natural metric $G$, and investigate the harmonicity of a vector field $V$ of $M$, thought as a map from $(M,g)$ to $(TM,G)$. We then apply this study to the Reeb vector field and, in particular, to Hopf vector fields on odd-dimensional spheres. 
\end{abstract}

\section {\red Introduction}

Let $(M,g)$ be an $n$-dimensional Riemannian manifold. Its tangent bundle $TM$, equipped with the so-called {\em Sasaki metric} $g^s$, has been extensively studied by several authors and in many different contexts.

In particular, given a compact Riemannian manifold $(M,g)$, Nouhaud \cite{N} considered the problem of determining harmonic sections of $(TM,g^s)$, that is, vector fields $V \in \mathfrak X(M)$ which define harmonic maps from $(M,g)$ to $(TM,g^s)$. She found the expression of the energy associated to $V$ and proved that parallel vector fields are all and the ones harmonic sections. Ishihara \cite{I} obtained independently the same result, giving also the explicit expression of the tension field associated to a vector field $V$.

Given a vector field $V$ over a compact Riemannian manifold, the energy associated to the map $V: (M,g) \rightarrow (TM,g^s)$ admits the following very simple expression \cite{N}, \cite{Wo}:
\begin{equation}\label{ener0}
E(V)=\frac{n}{2} {\rm vol} (M) + \frac{1}{2} \int _M ||\nabla V ||^2 dv_g,
\end{equation}
which, up to a constant, also corresponds to the {total bending} of $V$ \cite{W1}.
 
More recently, Gil-Medrano \cite{G1} proved that critical points of $E: \mathfrak X(M) \rightarrow \Real$, that is, the energy functional restricted to vector fields, are again parallel vector fields. Moreover, in the same paper she also determined the tension field associated to a unit vector field $V:(M,\bar g)\rightarrow(T_1 M,g^s)$, where $\bar g$ is a new Riemannian metric on $M$, and investigated the problem of determining when $V$ defines a harmonic map.

Investigating critical points of the energy associated to vector fields is an interesting purpose under different points of view. On the one hand, in many cases a distinguished vector field appears in a natural way, and it is worthwhile to see how the criticality of such a vector field is related to the geometry of the manifold. A well known example of this situation is given by the Reeb vector field $\xi$ of a contact metric manifold (\cite{P1}, \cite{P2}). On the other hand, vector fields determining harmonic maps, provide new and interesting examples of harmonic maps having as target some Riemannian manifolds endowed of an higly non-trivial geometry. For more details and the state of the art for criticality of vector fields, we can refer to the survey \cite{G2}.

The Sasaki metric $g^s$ has been the most investigated among all possible Riemannnian metrics on $TM$. However, in many different contexts such metrics showed a very "rigid"
behaviour. Moreover, $g^s$ represents only one possible choice inside a wide family of Riemannian metrics on $TM$, known as {\em Riemannian $g$-natural metrics}, which depend on several independent smooth functions from $\Real ^+$ to $\Real$. As their name suggests, those metrics arise from a very "natural"
construction starting from a Riemannian metric $g$ over $M$. The
introduction of $g$-natural metrics moves from the classification
of natural transformations of Riemannian metrics on manifolds to
metrics on tangent bundles \cite{KSe}, or equivalently, the
description of all first order natural operators $D:S_+^2T^*
\rightsquigarrow (S^2T^*)T$, transforming Riemannian metrics on
manifolds into metrics on their tangent bundles \cite{KoMSl} (see
also \cite{A}). Riemannian $g$-natural metrics have been
completely described in \cite{AS2}. They depend on six smooth
functions from $\Real ^+$ to $\Real$, special choices of which
give all the well known examples of Riemannian metrics on $TM$ as
$g^s$ itself, the Cheeger-Gromoll metric $g_{GC}$ and the metrics
investigated in \cite{O} (cf. Remark \ref{Exofmetrics}).

Both the rigidity of the Sasaki metric, and the fact mentioned above that several well known examples of Riemannian metrics on $TM$ are $g$-natural, make interesting to  investigate criticality of a vector field $V$, when $g^s$ is replaced by an arbitrary Riemannian $g$-natural metric $G$. In particular, the following questions arise:

\medskip
1) {\em  When $V:(M,g)\rightarrow(TM,G)$ defines a harmonic map?}  

\medskip
2) {\em When $V$ is a critical point for the energy $E$ restricted to vector fields?}  

\medskip\noindent
The aim of this paper is to answer the questions above. Note that in the study of Question 1, we shall find new examples of harmonic maps from $M$ to $TM$, defined by non-parallel vector fields (as Reeb vector fields and Hopf vector fields). The paper is organized in the following way. In Section 2, we shall recall the definition and basic properties of $g$-natural metrics on $TM$. The energy associated to $V:(M,g)\rightarrow(TM,G)$ when $M$ is compact, is explicitly calculated in Section 3, while in Section 4 we shall calculate the tension field associated to $V$. In Section 5, we shall determine some families of Riemannian $g$-natural metrics for which, as for $g^s$, parallel vector fields are all and the ones defining harmonic maps. In Section 6, we shall consider vector fields which are critical points for $E: \mathfrak X(M)\rightarrow \Real$, emphasizing the cases when this property is not equivalent to harmonicity of $V:(M,g)\rightarrow(TM,G)$. Finally, in Section 7 we shall apply our study to the case of the Reeb vector field $\xi$ of a contact metric manifold and, in particular, to Hopf vector fields on odd-dimensional spheres.

\section {\red Basic formulae on $g$-natural metrics on tangent bundles}

Let $(M,g)$ be an $n$-dimensional Riemannian manifold and $\nabla$ its Levi-Civita
connection. At any point $(x,u)$ of its {\em tangent bundle} $TM$, the tangent space of $TM$
splits into the horizontal and vertical subspaces with respect to
$\nabla$:
$$(TM)_{(x,u)}=\mathcal H _{(x,u)}\oplus \mathcal V _{(x,u)}.$$ \par

For any vector $X\in M_x$, there exists a unique vector $X^h  \in
\mathcal H _{(x,u)}$ (the \emph{horizontal lift} of $X$ to $(x,u)\in TM$), such that $\pi_*  X^h =X$, where $\pi:TM \rightarrow M$ is
the natural projection. The \emph{vertical lift} of a
vector $X\in M_x$ to $(x,u)\in TM$ is a vector $X^v  \in
\mathcal V_{(x,u)}$ such that $X^v (df) =Xf$, for all functions $f$ on $M$.
Here we consider $1$-forms $df$ on $M$ as functions on $TM$ (i.e.,
$(df)(x,u)=uf$). The map $X \to X^h$ is an isomorphism between the
vector spaces $M_x$ and $\mathcal H _{(x,u)}$. Similarly, the map $X \to
X^v$ is an isomorphism between $M_x$ and $\mathcal V _{(x,u)}$. Each tangent
vector $\tilde Z \in (TM)_{(x,u)}$ can be written in the form
$\tilde Z =X^h + Y^v$, where $X,Y \in M_x$ are uniquely determined
vectors. Horizontal and vertical lifts of vector fields on $M$ can be
defined in an obvious way and are uniquely defined vector fields
on $TM$. \par

We now write $F$ for the natural bundle with $FM =\pi^*(T^*\otimes T^*)M
\rightarrow M$. Then, we have $Ff(X_x,g_x)=(Tf.X_x,(T^* \otimes T^*)f.g_x)$
for all manifolds $M$, local diffeomorphisms $f$ of $M$, $X_x \in
T_x M$ and $g_x \in (T^*\otimes T^*)_x M$. The sections of the
canonical projection $FM \to M$ are called {\it $F$-metrics} in
literature. So, if we denote by $\oplus$ the fibered product of
fibered manifolds, then the $F$-metrics are mappings $TM \oplus TM
\oplus TM \to \Real$ which are linear in the second and the third
argument.\par

For a given $F$-metric $\delta$ on $M$, there are three
distinguished constructions of metrics on the tangent bundle $TM$
[KSe]:  \par

$(a)$ If $\delta$ is symmetric, then the \emph{Sasaki lift}
$\delta^s$ of $\delta$ is defined by
$$\left\{
\begin{array}{ll}
\delta^s_{(x,u)}(X^h,Y^h)=\delta(u;X,Y), &
\delta^s_{(x,u)}(X^h,Y^v)=0, \\
\delta^s_{(x,u)}(X^v,Y^h)=0, &
\delta^s_{(x,u)}(X^v,Y^v)=\delta(u;X,Y),
\end{array}   \right.$$
for all $X$, $Y$ $\in M_x$. When $\delta$ is non degenerate and
positive definite, so is $\delta^s$.\par

$(b)$ The \emph{horizontal lift} $\delta^h$ of $\delta$ is a
pseudo-Riemannian metric on $TM$, given by
$$\left\{
\begin{array}{ll}
\delta^h_{(x,u)}(X^h,Y^h)=0, &
\delta^h_{(x,u)}(X^h,Y^v)=\delta(u;X,Y),\\
\delta^h_{(x,u)}(X^v,Y^h)=\delta(u;X,Y),  &
\delta^h_{(x,u)}(X^v,Y^v)=0,
\end{array}    \right.$$
for all $X$, $Y$ $\in M_x$. If $\delta$ is positive definite, then
$\delta^s$ is of signature $(m,m)$. \par

$(c)$ The \emph{vertical lift} $\delta^v$ of $\delta$ is a
degenerate metric on $TM$, given by
$$\left\{
\begin{array}{ll}
\delta^v_{(x,u)}(X^h,Y^h)= \delta(u;X,Y),     &
\delta^v_{(x,u)}(X^h,Y^v)=0,     \\
\delta^v_{(x,u)}(X^v,Y^h)=0,            &
\delta^v_{(x,u)}(X^v,Y^v)=0,
\end{array}   \right.$$
for all $X$, $Y$ $\in M_x$. The rank of $\delta^v$ is exactly that
of $\delta$.  \\
If $\delta =g$ is a Riemannian metric on $M$, then these three
lifts of $\delta$ coincide with the three well-known classical
lifts of the metric $g$ to $TM$.\par

The three lifts above of \emph{natural} $F$-metrics generate the
class of $g$-natural metrics on $TM$. The introduction of $g$-natural metrics
moves from the description of all first order natural operators
$D:S_+^2T^* \rightsquigarrow (S^2T^*)T$, transforming Riemannian
metrics on manifolds into metrics on their tangent bundles, where
$S_+^2T^*$ and $S^2T^*$ denote the bundle functors of all
Riemannian metrics and all symmetric $(0,2)$-tensors over $n$-manifolds
respectively. For more details about the concept of naturality and related notions,
we can refer to \cite{KoMSl}.\par

Every section $G:TM \to (S^2T^*)TM$ is called a (possibly
degenerate) \emph{metric}. Then there is a bijective
correspondence between the triples of first order natural
$F$-metrics $(\zeta_1, \zeta_2, \zeta_3)$ and first order natural
(possibly degenerate) metrics $G$ on the tangent bundles given by
(cf. \cite{KSe}):
\begin{equation*}
G = \zeta_1^s +\zeta_2^h + \zeta_3^v.
\end{equation*} \par

Therefore, to find all first order natural operators $ S_+^2T^*
\rightsquigarrow (S^2T^*)T$ transforming Riemannian metrics on
manifolds into metrics on their tangent bundles, it suffices to
describe all first order natural $F$-metrics, i.e. first order
natural operators $ S_+^2T^* \rightsquigarrow (T,F)$. In this
sense, it is shown in \cite{KSe} (see also \cite{KoMSl}
and \cite{AS1}) that all first order natural $F$-metrics $\zeta$
in dimension $n>1$ form a family parametrized by two arbitrary
smooth functions $\alpha_0$, $\beta_0:\mathbb{R}^+ \to
\mathbb{R}$, where $\mathbb{R}^+$ denotes the set of all
nonnegative real numbers, in the following way: For every
Riemannian manifold $(M,g)$ and tangent vectors $u$, $X$, $Y \in
M_x$
\begin{equation}\label{exp-F-metr}
\zeta_{(M,g)}(u)(X,Y) =\alpha_0(g(u,u))g(X,Y)
+\beta_0(g(u,u))g(u,X)g(u,Y).
\end{equation}
If $n=1$, then the same assertion holds, but we can always choose
$\beta_0 =0$. In particular, all first order natural $F$-metrics
are symmetric.

We shall call a metric $G$ on $TM$, coming from $g$ by a first order natural operator $
S_+^2T^* \rightsquigarrow (S^2T^*)T$, a \emph{$g$-natural metric} \cite{AS2}. All $g$-natural metrics on the tangent bundle of a
Riemannian manifold $(M,g)$ are completely determined as follows:

\begin{Pro}[\cite{AS2}]\label{$g$-nat} 
Let $(M,g)$ be a Riemannian manifold and $G$ be a $g$-natural
metric on $TM$. Then there are six smooth functions $\alpha_i$,
$\beta_i:\mathbb{R}^+ \rightarrow \mathbb{R}$, $i=1,2,3$, such
that for every $u$, $X$, $Y\in M_x$, we have
\arraycolsep1.5pt
\begin{equation}\label{exp-$g$-nat}
 \left\lbrace
\begin{array}{rcl}
G_{(x,u)}(X^h,Y^h)& = & (\alpha_1+ \alpha_3)(r^2) g_x(X,Y)+ (\beta_1+ \beta_3)(r^2)g_x(X,u)g_x(Y,u),\\
G_{(x,u)}(X^h,Y^v)& = & \alpha_2 (r^2) g_x(X,Y)
                       +  \beta_2 (r^2) g_x(X,u)g_x(Y,u), \\
G_{(x,u)}(X^v,Y^h)& = & \alpha_2 (r^2) g_x(X,Y)
                       +  \beta_2 (r^2) g_x(X,u)g_x(Y,u), \\
G_{(x,u)}(X^v,Y^v)& = & \alpha_1 (r^2) g_x(X,Y)
                       + \beta_1 (r^2) g_x(X,u)g_x(Y,u),
\end{array}
\arraycolsep5pt \right.
\end{equation}
where $r^2 =g_x(u,u)$.\\ For $n=1$, the same holds with $\beta_i
=0$, $i=1,2,3$.
\end{Pro}

\begin{Nots}
In the sequel, we shall use the following notations:
\begin{itemize}
  \item $\phi_i(t) =\alpha_i(t) +t \beta_i(t)$,
  \item $\alpha(t) = \alpha_1(t) (\alpha_1+\alpha_3)(t) - \alpha_2
         ^2$(t),
  \item $\phi(t) = \phi_1(t) (\phi_1 +\phi_3)(t) -\phi_2^2(t)$,
\end{itemize}
for all $t \in \mathbb{R}^+$.
\end{Nots}

Riemannian $g$-natural metrics are characterized as follows:

\begin{Pro}[\cite{AS2}]\label{riem-nat} 
The necessary and sufficient conditions for a $g$-natural metric
$G$ on the tangent bundle of a Riemannian manifold $(M,g)$ to be
Riemann\-ian, are that the functions of Proposition \ref{$g$-nat},
defining $G$, satisfy the inequalities
\begin{equation}\label{cond-riem}
\left\lbrace
\begin{array}{ll}
\alpha_1(t) > 0, & \quad \phi_1(t) > 0,                      \\
\alpha(t) > 0, & \quad \phi(t)>0,
\end{array} \right.
\end{equation}
for all $t \in \mathbb{R}^+$.\par

For $n=1$, the system {\em (\ref{cond-riem})} reduces to $\alpha_1(t) >0$ and $\alpha(t)
>0$, for all $t \in \mathbb{R}^+$.
\end{Pro}

\begin{Conv}
a) In the sequel, when we consider an arbitrary Riemannian
$g$-natural metric $G$ on $TM$, we implicitly suppose that it is
defined by the functions $\alpha_i$, $\beta_i:\mathbb{R}^+
\rightarrow \mathbb{R}$, $i=1,2,3$, given in Proposition
\ref{$g$-nat} and satisfying (\ref{cond-riem}).\par

b) Unless otherwise stated, all real functions $\alpha_i$,
$\beta_i$, $\phi_i$, $\alpha$ and $\phi$ and their derivatives are
evaluated at $r ^2:=g_x(u,u)$.\par

c) We shall denote respectively by $R$ and $Q$ the {\em curvature tensor} and the {\em Ricci operator} of a Riemannian manifold $(M,g)$. The tensor $R$ is taken with the sign convention
$$R(X,Y)Z= \nabla _X \nabla _Y Z - \nabla _Y \nabla _X Z -\nabla _{[X,Y]}  Z,$$
for all vector fields $X,Y,Z$ on $M$.

\end{Conv}

\begin{Rk}\label{Exofmetrics}
In literature, there are some well known Riemannian metrics on the tangent bundle,  which turn out to be special cases of Riemannian $g$-natural metrics (satisfying (\ref{cond-riem})). In particular:
\begin{itemize}
\item the {\em Sasaki metric} $g^s$ is obtained for 
\begin{equation}\label{Sasa}
\alpha _1 (t)=1, \quad \alpha _2 (t)= \alpha _3 (t)= \beta _1 (t)=\beta _2 (t)= \beta _3 (t)=0.
\end{equation}
\item the {\em Cheeger-Gromoll metric} $g_{GC}$ \cite{CGr} is obtained when 
\begin{equation}\label{CheGro}
\alpha _2 (t)= \beta _2 (t)=0, \quad \alpha _1 (t)= \beta _1 (t)=-\beta _3 (t)= \frac{1}{1+t}, \quad \alpha _3 (t)= \frac{t}{1+t}.
\end{equation}

\item the two-parameters family of metrics investigated by Oproiu in \cite{O}, is obtained when there exist two smooth functions $v,w: \Real ^+ \rightarrow \Real$, such that (see \cite{AS2})
$$
\left\{
\begin{array}{l}
(\alpha _1 +\alpha _3)(t)=v(t/2), \quad (\beta _1 +\beta
_3)(t)=w(t/2), \\
\alpha _1 (t)=\frac{1}{v(t/2)}, \quad \beta
_1(t)=-\frac{w(t/2)}{v(t/2)[v(t/2)+tw(t/2)]}, \\
\alpha _2 (t)=\beta _2 (t)=0.
\end{array}
\right.
$$
\end{itemize}
Since $\alpha _2=\beta _2 =0 $, all these metrics are examples of Riemannian $g$-natural metrics on $TM$, for which horizontal and vertical distributions are mutually orthogonal.

\end{Rk}

The Levi-Civita connection $\bar \nabla$ of an arbitrary $g$-natural metric $G$ on $TM$, can be described as follows:

\begin{Pro}[\cite{AS1}]\label{lev-civ-con}
Let $(M,g)$ be a Riemannian manifold, $\nabla$ its Levi-Civita
connection and $R$ its curvature tensor. Let $G$ be a Riemannian
$g$-natural metric on $TM$. Then the Levi-Civita connection $\bar
\nabla$ of $(TM,G)$ is characterized by
$$ \arraycolsep1.5pt
\begin{array}{rcl}
(i) (\bar \nabla_{X^h}Y^h)_{(x,u)} & = &  (\nabla_X Y)^h _{(x,u)}
                   +h\{A(u;X_x,Y_x)\} + v\{B(u;X_x,Y_x)\}, \\
(ii) (\bar \nabla_{X^h}Y^v)_{(x,u)} & = &  (\nabla_X Y)^v _{(x,u)}
                   +h\{C(u;X_x,Y_x)\} + v\{D(u;X_x,Y_x)\}, \\
(iii) (\bar \nabla_{X^v}Y^h)_{(x,u)} & = &  h\{C(u;Y_x,X_x)\} +
                      v\{D(u;Y_x,X_x)\}, \\
(iv) (\bar \nabla_{X^v}Y^v)_{(x,u)} & = &
                      h\{E(u;X_x,Y_x)\} + v\{F(u;X_x,Y_x)\},
\end{array}
\arraycolsep5pt $$ for all vector fields $X$, $Y$ on $M$ and
$(x,u) \in TM$. Here, $h\{ \cdot \}$ and $v\{ \cdot \}$ respectively denote the horizontal and vertical lifts of a vector tangent to $M$ and, for all $x \in M$ and vectors $u$, $X_x$, $Y_x$ tangent to $M$ at $x$, $A$, $B$, $C$, $D$, $E$ and $F$ are defined as follows: 
$$
\begin{array}{lcl}
A(u;X_x,Y_x)  & = & A_1 [R_x(X_x,u)Y_x +R_x(Y_x,u)X_x]+A_2 [g_x(Y_x,u)X_x + g_x(X_x,u)Y_x] \vphantom{\displaystyle\frac{a}{a}} \\ & &+A_3 g_x(R_x(X_x,u)Y_x,u)u + A_4  g_x(X_x,Y_x)u + A_5 g_x(X_x,u)g_x(Y_x,u)u,  \vphantom{\displaystyle\frac{a}{a}} \end{array}
\arraycolsep5pt$$
where
\begin{equation}\label{Ai}
\begin{array}{lcl}
\vphantom{\displaystyle\frac{A}{A}} A_1 &=& -\frac{\alpha_1 \alpha_2}{2\alpha} , \\
\vphantom{\displaystyle\frac{A}{A}} A_2 &=& \frac{\alpha_2 (\beta_1 +\beta_3)}{2\alpha}, \\ \vphantom{\displaystyle\frac{A}{A}} A_3 &=& \frac{ \alpha_2\{ \alpha_1[\phi_1 (\beta_1 +\beta_3)
                 -\phi_2 \beta_2] +\alpha_2(\beta_1 \alpha_2 - \beta_2 \alpha_1)\}}{\alpha \phi} , \\
\vphantom{\displaystyle\frac{A}{A}} A_4 &=& \frac{ \phi_2 (\alpha_1 + \alpha_3)^\prime }{\phi} , \\
\vphantom{\displaystyle\frac{A}{A}} A_5 &=& \frac{\alpha \phi_2(\beta_1 +\beta_3)^\prime + (\beta_1 +\beta_3)\{ \alpha_2 [\phi_2 \beta_2 - \phi_1 (\beta_1 +\beta_3)] + (\alpha_1 +\alpha_3)(\alpha_1 \beta_2
                 - \alpha_2 \beta_1)\}}{\alpha \phi} , \vphantom{\displaystyle\frac{A}{A}}
\end{array}
\arraycolsep5pt
\end{equation}

$$\arraycolsep1.5pt
\begin{array}{lcl}
B(u;X_x,Y_x) & = &  B_1 R_x(X_x,u)Y_x +B_2 R_x(X_x,Y_x)u +B_3 [g_x(Y_x,u)X_x +g_x(X_x,u)Y_x] \vphantom{\displaystyle\frac{a}{a}} \\
         & & + B_4 g_x(R_x(X_x,u)Y_x,u)u + B_5 g_x(X_x,Y_x)u +B_6 g_x(X_x,u)g_x(Y_x,u)u , \vphantom{\displaystyle\frac{a}{a}} 
         \end{array}
\arraycolsep5pt$$
where
\begin{equation}\label{Bi}
\arraycolsep1.5pt
\begin{array}{lcl}
\vphantom{\displaystyle\frac{A}{A}} B_1 &=&  \frac{\alpha_2 ^2}{\alpha}, \\
\vphantom{\displaystyle\frac{A}{A}} B_2 &=& -\frac{\alpha_1 (\alpha_1 +\alpha_3)}{2\alpha}, \\
\vphantom{\displaystyle\frac{A}{A}} B_3 &=& - \frac{(\alpha_1 +\alpha_3) (\beta_1 +\beta_3)}
               {2\alpha}, \\
\vphantom{\displaystyle\frac{A}{A}} B_4 &=&  \frac{\alpha_2 \{ \alpha_2[\phi_2 \beta_2- \phi_1 (\beta_1
               +\beta_3)] +(\alpha_1 +\alpha_3)(\beta_2
               \alpha_1  - \beta_1 \alpha_2)\}}{\alpha\phi}, \\
\vphantom{\displaystyle\frac{A}{A}} B_5 &=& - \frac{(\phi_1
               +\phi_3)(\alpha_1+ \alpha_3)^\prime}{\phi}, \\    
\vphantom{\displaystyle\frac{A}{A}} B_6 &=& \frac{-\alpha (\phi_1 +\phi_3)(\beta_1 +\beta_3)^\prime +(\beta_1 +\beta_3) \{ (\alpha_1 +\alpha_3)[(\phi_1 +\phi_3)
                \beta_1- \phi_2 \beta_2] 
                + \alpha_2[\alpha_2 (\beta_1 +\beta_3)
                -(\alpha_1 +\alpha_3) \beta_2 ] \}}{\alpha\phi} , 
\end{array}
\arraycolsep5pt
\end{equation}

$$\arraycolsep1.5pt
\begin{array}{lcl}
C(u;X_x,Y_x) & = &  C_1 R(Y_x,u)X_x +C_2 g_x(X_x,u)Y_x +C_3 g_x(Y_x,u)X_x  \vphantom{\displaystyle\frac{A}{A}} \vphantom{\displaystyle\frac{A}{A}} \\& &+C_4 g_x(R_x(X_x,u)Y_x,u)u +C_5 g_x(X_x,Y_x)u +C_6 g_x(X_x,u)g_x(Y_x,u) u , 
\end{array}
\arraycolsep5pt$$

where
\begin{equation}\label{Ci}
\arraycolsep1.5pt
\begin{array}{lcl}
\vphantom{\displaystyle\frac{A}{A}} C_1 &=&  -\frac{\alpha_1^2}{2 \alpha}, \\
\vphantom{\displaystyle\frac{A}{A}} C_2 &=& -\frac{\alpha_1 (\beta_1 +\beta_3)}{2 \alpha}, \\
\vphantom{\displaystyle\frac{A}{A}} C_3 &=& \frac{\alpha_1 (\alpha_1 +\alpha_3)^\prime
               -\alpha_2(\alpha_2 ^\prime -\frac{\beta_2}{2})}{\alpha}, \\
\vphantom{\displaystyle\frac{A}{A}} C_4 &=&  \frac{ \alpha_1 \{ \alpha_2 (\alpha_2 \beta_1 -
               \alpha_1 \beta_2) +\alpha_1 [\phi_1 (\beta_1
               +\beta_3) -\phi_2 \beta_2]\}}{2\alpha \phi}
              , \\
\vphantom{\displaystyle\frac{A}{A}} C_5 &=&   \frac{\phi_1 (\beta_1 +\beta_3)
               +\phi_2(2 \alpha_2 ^\prime -\beta_2)}{2\phi} , \\    
\vphantom{\displaystyle\frac{A}{A}} C_6 &=&
\frac{\alpha \phi_1  (\beta_1 +\beta_3)^\prime
+ \{ \alpha_2(\alpha_1 \beta_2 -\alpha_2 \beta_1)
+\alpha_1[\phi_2 \beta_2 -(\beta_1 +\beta_3) \phi_1] \}
              [(\alpha_1 +\alpha_3)^\prime
              +\frac{\beta_1 +\beta_3}2] }{\alpha \phi} \\
 & & + \frac{\{\alpha_2 [\beta_1 (\phi_1 +\phi_3) -\beta_2 \phi_2]
 -\alpha_1 [\beta_2 (\alpha_1 +\alpha_3)
 -\alpha_2(\beta_1 +\beta_3)]\}(\alpha_2 ^\prime
 -\frac{\beta_2}{2})}{\alpha \phi} 
\end{array}
\arraycolsep5pt
\end{equation}

$$\arraycolsep1.5pt
\begin{array}{lcl}
D(u;X_x,Y_x) & = & D_1 R_x(Y_x,u)X_x +D_2 g_x(X_x,u)Y_x + D_3 g_x(Y_x,u)X_x  \vphantom{\displaystyle\frac{A}{A}} \\
         &  &+D_4 g_x(R_x(X_x,u)Y_x,u)u+ D_5 g_x(X_x,Y_x)u +D_6 g_x(X_x,u)g_x(Y_x,u) u ,
\end{array}
\arraycolsep5pt$$

where
\begin{equation}\label{Di}
\arraycolsep1.5pt
\begin{array}{lcl}
\vphantom{\displaystyle\frac{A}{A}} D_1 &=& \frac{\alpha_1 \alpha_2}{2 \alpha}, \\
\vphantom{\displaystyle\frac{A}{A}} D_2 &=& \frac{\alpha_2 (\beta_1 +\beta_3)}{2 \alpha}, \\
\vphantom{\displaystyle\frac{A}{A}} D_3 &=& \frac{ -\alpha_2 (\alpha_1 +\alpha_3)^\prime  +(\alpha_1 +\alpha_3) (\alpha_2 ^\prime -\frac{\beta_2}{2})}{\alpha}, \\
\vphantom{\displaystyle\frac{A}{A}} D_4 &=&  \frac{\alpha_1 \{ (\alpha_1 +\alpha_3)(\alpha_1 \beta_2 -
              \alpha_2 \beta_1) +\alpha_2 [\phi_2 \beta_2 -\phi_1 (\beta_1
               +\beta_3)]\}}{2\alpha \phi} \\
      \vphantom{\displaystyle\frac{A}{A}} D_5 &=&   -\frac{\phi_2 (\beta_1 +\beta_3)
               +(\phi_1 +\phi_3)(2\alpha_2 ^\prime -\beta_2)}{2\phi} , \\    
\vphantom{\displaystyle\frac{A}{A}} D_6 &=&  \frac{-\alpha \phi_2 (\beta_1 +\beta_3)^\prime
               +\{ (\alpha_1 +\alpha_3)(\alpha_2 \beta_1 -\alpha_1 \beta_2)
                 +\alpha_2[\phi_1(\beta_1 +\beta_3) -\phi_2 \beta_2]\}
               [(\alpha_1 +\alpha_3)^\prime +\frac{\beta_1 +\beta_3}2]}{\alpha \phi} \\
 & & + \frac{\{ (\alpha_1 +\alpha_3)[\beta_2 \phi_2 -\beta_1 (\phi_1
               +\phi_3)]+\alpha_2 [\beta_2 (\alpha_1 +\alpha_3)-\alpha_2(\beta_1 +\beta_3)]\}(\alpha_2 ^\prime  -\frac{\beta_2}{2})}{\alpha \phi}  
\end{array}
\arraycolsep5pt
\end{equation}

$$\arraycolsep1.5pt
\begin{array}{lcl}
E(u;X_x,Y_x) & = &  E_1 [g_x(Y_x,u)X_x +g_x(X_x,u)Y_x] + E_2 g_x(X_x,Y_x)u+E_3 g_x(X_x,u)g_x(Y_x,u) u  
\end{array}
\arraycolsep5pt$$

where
\begin{equation}\label{Ei}
\arraycolsep1.5pt
\begin{array}{lcl}
\vphantom{\displaystyle\frac{A}{A}} E_1 &=& \frac{\alpha_1 (\alpha_2^\prime
               +\frac{\beta_2}{2}) -\alpha_2 \alpha_1^\prime}{\alpha} , \\
\vphantom{\displaystyle\frac{A}{A}} E_2 &=& \frac{\phi_1 \beta_2 -\phi_2(\beta_1 -\alpha_1 ^\prime)}{\phi} , \\
\vphantom{\displaystyle\frac{A}{A}} E_3 &=& \frac{\alpha (2\phi_1 \beta_2^\prime
               -\phi_2 \beta_1^\prime) +2\alpha_1^\prime
               \{ \alpha_1[\alpha_2 (\beta_1 +\beta_3) - \beta_2 (\alpha_1 +\alpha_3)] +\alpha_2[\beta_1(\phi_1
               +\phi_3) -\beta_2 \phi_2]\} }{\alpha \phi} \\
                & & + \frac{(2\alpha_2 ^\prime +\beta_2) \{ \alpha_1[\phi_2 \beta_2
               -\phi_1(\beta_1 +\beta_3)] +\alpha_2 (\alpha_1 \beta_2
                -\alpha_2 \beta_1)\}}{\alpha \phi}  
\end{array}
\arraycolsep5pt
\end{equation}

$$\arraycolsep1.5pt
\begin{array}{lcl}
F(u;X_x,Y_x) & = &  F_1 [g_x(Y_x,u)X_x +g_x(X_x,u)Y_x] +F_2 g_x(X_x,Y_x)u +F_3 g_x(X_x,u)g_x(Y_x,u) u 
\end{array}
\arraycolsep5pt$$

where
\begin{equation}\label{Fi}
\arraycolsep1.5pt
\begin{array}{lcl}
\vphantom{\displaystyle\frac{A}{A}} F_1 &=&  \frac{-\alpha_2 (\alpha_2^\prime
               +\frac{\beta_2}{2}) +(\alpha_1 +\alpha_3) \alpha_1^\prime}{\alpha}, \\
\vphantom{\displaystyle\frac{A}{A}} F_2 &=& \frac{ (\phi_1 +\phi_3)(\beta_1 -\alpha_1 ^\prime)
               -\phi_2 \beta_2}{\phi} \\
\vphantom{\displaystyle\frac{A}{A}} F_3 &=& \frac{\alpha [(\phi_1 +\phi_3) \beta_1^\prime
               -2\phi_2 \beta_2^\prime] +2\alpha_1^\prime \{ \alpha_2
               [\beta_2 (\alpha_1 +\alpha_3)- \alpha_2 (\beta_1 +\beta_3)] + (\alpha_1 +\alpha_3)
               [\beta_2 \phi_2 -\beta_1(\phi_1+\phi_3)]\}}{\alpha \phi} \\
                & & + \frac{(2\alpha_2 ^\prime +\beta_2) \{ \alpha_2[\phi_1(\beta_1
               +\beta_3) -\phi_2 \beta_2] +(\alpha_1 +\alpha_3)(\alpha_2\beta_1
               -\alpha_1 \beta_2)\}}{\alpha \phi}  
\end{array}
\arraycolsep5pt
\end{equation}

For $n=1$, the same holds with $\beta_i =0$, $i=1,2,3$.
\end{Pro}

\section {\red The energy of a vector field $V: (M,g) \rightarrow (TM, G)$} \label{energy}

We shall first discuss geometric properties of the map $V: (M,g) \rightarrow (TM,G)$ defined by a vector field $V \in \mathfrak X(M)$.
It is well known that if $TM$ is equipped with the Sasaki metric $g^s$, then $V$ defines an isometry $V:(M,g) \rightarrow (TM,g^s)$, that is, it satisfies $V^* g^s=g$, if and only if $V$ is parallel. We now replace $g^s$ by an arbitrary Riemannian $g$-natural metric $G$. Since $V_* X=X^h + (\nabla _X V)^v$ for any vector field $X$, from (\ref{exp-$g$-nat}) we obtain
\begin{eqnarray}\label{isom}
 (V^* G)(X,Y) &=& G(X^h+ (\nabla _X V)^v , Y^h + (\nabla _Y V)^v ) \nonumber \\ &=& (\alpha _1 +\alpha _3) (r^2) g(X,Y) + (\beta _1 +\beta _3 )(r^2) g(X,V)g(Y,V) 
 \\ & &  + \alpha _2 (r^2) \left[ g(X,\nabla_Y V)+g(Y,\nabla _X V) \right]  \nonumber  \\
& & + \beta _2 (r^2)  \left[ g(X,V)g(\nabla_Y V,V)+g(Y,V)g(\nabla_X V,V) \right] \nonumber \\& & + \alpha _1 (r^2) g(\nabla _X V, \nabla _Y V) +
 \beta _1 (r^2) g(\nabla _X V,V)g(\nabla _Y V,V), \nonumber 
\end{eqnarray}
for all vector fields $X,Y$, where $r=||V||$ is a smooth function from $M$ to $\Real ^+$. Note that by (\ref{isom}), in general $V^* G$ also depends on the length of $V$. 

In particular, under the assumption $\beta _1 + \beta _3 =0$, (determining a very large family of $g$-natural metrics, which includes $g^s$ and depends on five smooth functions $\alpha _1$, $\alpha _2$, $\alpha _3$, $\beta _1$  and $\beta _2$), from (\ref{isom}) we easily get the following   

\begin{Pro}\label{teor1}
Let $G$ be a Riemannian $g$-natural metric and $V \in \mathfrak X(M)$. 

\medskip
1) If $\beta _1 +\beta _3 =0$, then 

\medskip\noindent
a) $\nabla V =0$ implies that $V^* G$ is homothetic to $g$, with homotethy factor 
$(\alpha _1 +\alpha _3) (\rho)$, where $\rho =||V||^2$ is constant. In particular, $V$ is an isometry when in addition $(\alpha _1 +\alpha _3) (\rho)=1$.

\medskip\noindent
b) If $M$ is compact and $V$ has constant length $||V||=\sqrt{\rho}$, then

$$V^* G =(\alpha _1 +\alpha _3) (\rho) \, g \Leftrightarrow \nabla V=0. $$ 

\medskip
2) If $\alpha _2 =\beta _1 =\beta _2 = \beta _3 =0$, then $V$ is parallel if and only if $V^* G =(\alpha _1 +\alpha _3) (r^2) \, g$.

\end{Pro}

\medskip\noindent
{\bf Proof.} {\em 1)}: {\em a)}: it  follows at once from (\ref{isom}), rewritten when $\beta _1 +\beta _3 =0$ and $\nabla V=0$. 

{\em b)}: the "if" part follows from {\em a)}. For the "only if" part, we consider a local orthonormal basis $\{e_i\}$ on $M$ and apply (\ref{isom}) to pairs $(e_i,e_i)$, for all $i$=1,..,n. Taking into account the fact that $V^* G =(\alpha _1 +\alpha _3) (\rho) \, g$ and summing up over $i$, we easily get
\begin{equation}\label{temp}
2 \alpha _2 (\rho) {\rm div} V +  \alpha _1 (\rho) || \nabla V ||^2 =0.
\end{equation}
Since $M$ is compact, we can integrate (\ref{temp}) over $M$ and we obtain
\begin{equation}\label{temp1}
\alpha _1 (\rho) \int _M ||\nabla V||^2 dv_g =0,
\end{equation}
because $\rho$ is constant and $\int _M {\rm div} V dv_g =0$. By (\ref{cond-riem}), $\alpha _1 >0$. So, (\ref{temp1}) yields $||\nabla V ||^2 =0$, that is, $V$ is parallel.

\medskip\noindent
{\em 2)}: the "if" part follows directly from {\em a)}. For the "only if" part, it is enough to rewrite (\ref{isom}) 
for $\alpha _2 =\beta _1 =\beta _2 = \beta _3 =0$ and $V^* G =(\alpha _1 +\alpha _3) (r^2) \, g$, and we get
$$\alpha _1 (r^2) g(\nabla _X V, \nabla _Y V) =0$$ for all vector fields $X,Y$. Since $\alpha _1 >0$, we then have $\nabla V=0$ $\Box$

\bigskip\noindent
In order to provide some examples, note that, by (\ref{Sasa}) and (\ref{CheGro}), the Sasaki metric $g^s$ and the Cheeger -Gromoll metric $g_{CG}$ on $TM$ satisfy conditions listed at points {\em 2)} and {\em 1)} of Proposition \ref{teor1}, respectively.

\medskip
Next, let $f: (M,g) \rightarrow (M',g')$ be a smooth map between Riemannian manifolds, with $M$ compact. The  {\em energy} of $f$ is defined as the integral 
$$ E(f) := \int _M e(f) dv_g $$
where $e(f)= \frac{1}{2} ||f_*||^2= \frac{1}{2} {\rm tr}_g f^* g'$ is  the so-called {\em energy density} of $f$. With respect to a local orthonormal basis of vector fields $\{e_1,..,e_n \}$ on $M$, it is possible to express the energy density as $e(f)= \frac{1}{2} \sum _{i=1} ^n g'(f_*e_i,f_*e_i)$. Critical points of the energy functional on $C^{\infty} (M,M')$ are known as {\em harmonic maps}. They have been characterized in \cite{ES} as maps having vanishing {\em tension field} $\tau (f)={\rm tr} \nabla df$. When $(M,g)$ is a general Riemannian manifold (including the non-compact case), a map $f: (M,g) \rightarrow (M',g')$ is said to be harmonic if $\tau (f)=0$. For further details about the energy functional, we can refer to \cite{EL1},\cite{U}. 

Let now $(M,g)$ be a compact Riemannian manifold of dimension $n$ and $(TM,G)$ its tangent bundle, equipped with an arbitrary Riemannian $g$-natural metric $G$. Each vector field $V \in \mathfrak X (M)$ defines a smooth map $V: (M,g) \rightarrow (TM,G)$, $p \mapsto V_p $. By definition, the {\em energy} $E(V)$ of $V$ is the energy associated to the corresponding map $V: (M,g) \rightarrow (TM, G)$. Therefore, $E(V)=\int _M e(V) dv_g$, where the density function $e(V)$ is given by 
\begin{equation}
	e _p (V) = \frac{1}{2} || V_{*p}|| ^2 = \frac{1}{2} {\rm tr}_g  (V^{*} \, G)_p = \frac{1}{2} \sum _{i=1} ^n (V^{*} G)_p (e_i,e_i), 
\end{equation}
$\{ e_1,..,e_n \}$ being any local orthonormal basis of vector fields defined in a neighborhood of $p$.
 Using formulae (i)-(iv) of Proposition \ref{lev-civ-con},  we then have
\begin{eqnarray*}
e (V) &=& \frac{1}{2} \sum _{i=1} ^n G _{V} (V_{*} \, e_i, V_{*} \, e_i)=
 \frac{1}{2} \sum _{i=1} ^n G _{V} (e_i ^h +(\nabla _{e_i} V )^v, e_i ^h +(\nabla _{e_i} V )^v) \\
 & =& \frac{1}{2} \sum _{i=1} ^n \left\{(\alpha _1 +\alpha _3 )(r^2) g(e_i,e_i)+
 (\beta _1 +\beta _3 )(r^2) g(e_i,V) ^2 \right. \nonumber \\
 & &\qquad \quad +2\alpha _2 (r^2) g(e_i,\nabla _{e_i} V )+ 2\beta _2 (r^2) g(e_i,V) g(\nabla _{e_i} V ,V) \nonumber  
 \\ & & \qquad \quad \left. + \alpha _1 (r^2) g(\nabla _{e_i} V ,\nabla _{e_i} V )+
 \beta _1 (r^2) g(\nabla _{e_i} V ,V) ^2 \right\} \vphantom{\displaystyle\frac{A}{A}}
 \end{eqnarray*}
where $r = ||V||$ and so,
\begin{eqnarray}\label{dens}
\qquad e (V) &=& \frac{1}{2} \left\{ \vphantom{\displaystyle \frac{A}{A}} n(\alpha _1 +\alpha _3 )(r^2) + (\beta _1 +\beta _3 )(r^2) r^2 + 2\alpha _2 (r^2) {\rm div} (V) \right.   \vphantom{\displaystyle\frac{a}{a}} \\
 & & \left. + 2\beta _2 (r^2) V (r^2) + \alpha _1 (r^2) ||\nabla V||  ^2 +\frac{1}{4}
 \beta _1 (r^2) || {\rm grad} \, r ^2 ||^2 \right\} \nonumber \vphantom{\displaystyle\frac{a}{a}}.
 \end{eqnarray}
We now assume that $M$ is compact and we rewrite $ E(V) = \int _M e  (V) dv_g$ for some special kinds of vector fields. More precisely, we consider vector fields of constant length and, as a special case, parallel vector fields.

For any constant $\rho >0$, we put
$$\mathfrak X ^{\rho}(M) =\{V \in \mathfrak X (M): ||V||^2 = \rho \}.$$
So, if $V \in \mathfrak X ^{\rho}(M)$, then $V$ has constant length satisfying $||V||^2=\rho$. By (\ref{dens}) and taking into account the definition of $\phi _i$ given in Notations 1, we easily get that the energy of $V$ is given by
\begin{equation}\label{densunit}
E (V) = \frac{1}{2} [ (n-1)(\alpha _1 +\alpha _3 )+ 
 \phi _1 +\phi _3 ] (\rho) \cdot {\rm vol}(M,g) + \frac{1}{2} \alpha_1 (\rho) \cdot \int _M ||\nabla V||^2 dv_g. 
 \end{equation}
Since $\alpha _1 >0$, (\ref{densunit}) implies that 

\begin{equation}\label{inequnit}
E (V) \geq \frac{1}{2} [ (n-1)(\alpha _1 +\alpha _3 )+ 
 \phi _1 +\phi _3 ] (\rho) \cdot {\rm vol}(M,g) >0,
 \end{equation}
for all $V \in \mathfrak X ^{\rho}(M)$. (The last inequality follows at once from Notations 1 and  (\ref{cond-riem})). The equality holds in (\ref{inequnit}) if and only if $V$ is  parallel. Therefore, we have the following

\begin{Th}\label{teor2}
Let $(M,g)$ be a compact Riemannian manifold. Equipping $TM$ with an arbitrary Riemannian $g$-natural metric $G$, a vector field $V \in \mathfrak X ^{\rho}(M)$ is an absolute minimum for the energy $E: \mathfrak X ^{\rho}(M) \rightarrow \Real$ restricted to $\mathfrak X ^{\rho}(M)$ if and only if $V$ is parallel.
\end{Th}

\noindent
In particular, from Proposition \ref{teor1} and Theorem \ref{teor2} it follows 

\begin{Cor}\label{cor1}
Let $(M,g)$ be a compact Riemannian manifold and $V \in \mathfrak X ^{\rho}(M)$. With respect to a Riemannian $g$-natural metric $G$ satisfying $\beta _1 +\beta _3 =0$, the following assertions are equivalent:

\medskip
(i) $V$ is an absolute minimum for the energy $E: \mathfrak X ^{\rho}(M) \rightarrow \Real$,

\medskip
(ii) $V$ is parallel,

\medskip
(iii) $V^* g = (\alpha _1 +\alpha _3) (\rho) \, g$ (that is, $V:(M,g) \rightarrow (TM,G)$ is a homothetic immersion).
\end{Cor}

\noindent
It is worth mentioning that Corollary \ref{cor1} applies to both the Sasaki metric $g^s$ and the Cheeger -Gromoll metric $g_{CG}$ of $TM$.

Note that a parallel vector field $V$ necessarily has constant length. In fact, for all $X \in \mathfrak X (M)$ we have $2X (||V||^2) =g (\nabla _{X} V,V) =0$ . When $V$ is parallel, from (\ref{densunit}) (or (\ref{dens})) we have
\begin{equation}\label{denspar}
E(V) = \frac{1}{2} [ (n-1)(\alpha _1 +\alpha _3 )+ 
 \phi _1 +\phi _3 ] (\rho) \cdot {\rm vol}(M,g),
 \end{equation}
where $||V||^2=\rho$. By (\ref{denspar}), a parallel vector field $V$ 
is a critical point for the energy restricted to the set  
$$\mathfrak X _{\mathcal P} (M) =\{V \in \mathfrak X (M): \nabla V =0 \}$$
of all {\em parallel} vector fields, if and only if   
\begin{equation}\label{condpar}
[ (n-1)(\alpha _1 +\alpha _3 )+  \phi _1 +\phi _3 ]' (\rho)=0 , 
 \end{equation}
that is, $\rho =||V||^2$ is a critical point of the function $[ (n-1)(\alpha _1 +\alpha _3 )+ 
 \phi _1 +\phi _3 ]$. As we shall see in the next Section, (\ref{condpar}) is also a sufficient condition for a parallel vector field $V$ to define a harmonic map $V: (M,g) \rightarrow (TM, G)$.

\section {\red The tension field associated to $V: (M,g) \rightarrow (TM, G)$} \label{tension}

Let $(M,g)$ be a Riemannian manifold and $V \in \mathfrak X (M)$. The {\em tension field} associated to the map $V : (M,g) \rightarrow (TM,G)$, is defined as
\begin{equation}\label{deftens}
\arraycolsep2pt
\begin{array}{rcl}
\tau (V) : M & \rightarrow & V^{-1} (TTM) \\ \\
p & \mapsto & {\rm tr} (\nabla dV)_p .
\end{array}
 \end{equation}
Let $p$ be a point of $M$ and $\{e_1,..,e_n\}$ a local orthonormal basis of vector fields, defined in a neighborhood of $p$. By (\ref{deftens}), we have
\begin{eqnarray}\label{tens0}
\tau _p (V) &=& \sum _{i=1} ^n (\bar \nabla dv) (e_i,e_i)(p) =  \sum _{i=1} ^n \{ \bar \nabla _{V_{*} e_i} V_{*} e_i - V_{*} (\nabla _{e_i} e_i)\}(p)
\\
&=& \sum _{i=1} ^n \left\{ \bar \nabla _{e_i ^h + (\nabla _{e_i} V)^v} (e_i ^h + (\nabla _{e_i V})^v) -  (\nabla _{e_i} e_i) ^v - (\nabla _{\nabla _{e_i} e_i} V) ^v \right\}(p) \nonumber
\\
&=& \sum _{i=1} ^n \left\{ \bar \nabla _{e_i ^h} e_i ^h + \bar \nabla _{e_i ^h} (\nabla _{e_i} V)^v + \bar \nabla _{(\nabla _{e_i} V)^v} e_i ^h + \bar \nabla _{(\nabla _{e_i} V)^v} (\nabla _{e_i} V)^v \right. \nonumber \\ & & \left. \qquad -  (\nabla _{e_i} e_i) ^v - (\nabla _{\nabla _{e_i} e_i} V) ^v \vphantom{\displaystyle \nabla _{(\nabla _{e_i} V)^v}} \right\}(p). \nonumber
\end{eqnarray}
Hence, taking into account formulae of Proposition \ref{lev-civ-con} for the Levi-Civita connection of an arbitrary Riemannian $g$-natural metric $G$ on $TM$, from (\ref{tens0}) we easily get
\begin{eqnarray}\label{tfield}
\tau _p (V)  & = &  \left\{ \vphantom{\displaystyle\frac{A}{A}} -2 A_1 QV +2 C_1 {\rm tr} [R(\nabla _{\cdot} V,V)\cdot ]
+C_3 \sum _{i=1} ^n  e_i (r^2) e_i  \right.
\\ & &+2C_2 \nabla _V V +E_1 \sum _{i=1} ^n  e_i (r^2) \nabla _{e_i} V +\left[ \vphantom{\displaystyle\frac{A}{A}} 2A_2 -A_3 g(QV,V)+nA_4  \right.\nonumber \\ & & +A_5 r^2  +2C_4 g({\rm tr} [R(\nabla _{\cdot} V,V)\cdot ],V) +2C_5 {\rm div}V
  + C_6 V (r^2)  \nonumber \\ & & \left. \left. +E_2 ||\nabla V||^2 + \frac{1}{4} E_3 \sum _{i=1} ^n  [e_i (r^2)]^2 \right] V\right\} ^h _p \nonumber \\
  & & +  \left\{ \vphantom{\displaystyle\frac{A}{A}} -\bar \Delta V - B_1 QV +2 D_1 {\rm tr} [R(\nabla _{\cdot} V,V)\cdot ]
+D_3 \sum _{i=1} ^n  e_i (r^2) e_i \right. \nonumber
\\ & &+2D_2 \nabla _V V +F_1 \sum _{i=1} ^n  e_i (r^2) \nabla _{e_i} V +\left[ \vphantom{\displaystyle\frac{A}{A}} 2B_3 -B_4 g(QV,V)+nB_5 \right. \nonumber \\ & & +B_6 r^2 +2D_4 g({\rm tr} [R(\nabla _{\cdot} V,V)\cdot ],V) +2D_5 {\rm div}V  + D_6 V (r^2) \nonumber\\ & & \left.\left. 
  +F_2 ||\nabla V||^2 + \frac{1}{4} F_3 \sum _{i=1} ^n  [e_i (r^2)]^2 \right] V\right\} ^v _p , \nonumber
 \end{eqnarray}
where $r= ||V||$, $A _i,...F _i$ are evaluated at $r^2$ and $\bar \Delta V = -{\rm tr} \nabla ^2 V = -\sum _i \left(\nabla _{e_i} \nabla _{e_i} V - \nabla _{\nabla _{e_i} e_i} V \right)$ is the socalled {\em rough Laplacian} of $(M,g)$ calculated at $V$. Therefore, for the smooth map $V:(M,g) \rightarrow (TM,G)$ defined by a vector field $V \in \mathfrak X (M)$, by 
(\ref{tfield}) we obtain the following

\begin{Th}\label{harmtens}
Let $(M,g)$ be a compact Riemannian manifold. A vector field $V \in \mathfrak X (M)$  defines a harmonic map $V:(M,g) \rightarrow (TM,G)$ if and only if 
\begin{eqnarray}\label{htension}
& & \tau _h (V) = -2 A_1 QV +2 C_1 {\rm tr} [R(\nabla _{\cdot} V,V)\cdot ]
+C_3 {\rm grad} \, r ^2 + E_1 \nabla _{{\rm grad} \, r ^2} V 
\\ & & +2C_2 \nabla _V V +\left[ \vphantom{\displaystyle\frac{A}{A}} 2A_2 -A_3 g(QV,V)+nA_4  +A_5 r^2  +2C_4 g({\rm tr} [R(\nabla _{\cdot} V,V)\cdot ],V) 
    \nonumber \right. \\ & & \left. +2C_5 {\rm div}V + C_6 V (r^2)+E_2 ||\nabla V||^2 + \frac{1}{4}E_3 \left| \left| {\rm grad} \, r ^2 \right| \right|^2 \right] V =0 \nonumber 
\end{eqnarray}
and 
\begin{eqnarray}\label{vtension}
  & & \tau _v (V) =  -\bar \Delta V - B_1 QV +2 D_1 {\rm tr} [R(\nabla _{\cdot} V,V)\cdot ]
+D_3 {\rm grad} \, r ^2 +F_1 \nabla _{{\rm grad} \, r ^2} V 
\\ & &+2D_2 \nabla _V V +\left[ \vphantom{\displaystyle\frac{A}{A}} 2B_3 -B_4 g(QV,V)+nB_5 +B_6 r^2 +2D_4 g({\rm tr} [R(\nabla _{\cdot} V,V)\cdot ],V)  \nonumber \right. \\ & & \left. +2D_5 {\rm div}V   + D_6 V (r^2)
  +F_2 ||\nabla V||^2 + \frac{1}{4} F_3 \left| \left| {\rm grad} \, r ^2 \right| \right|^2 \right] V=0, \nonumber
\end{eqnarray}
where, for all points $p \in M$, $\tau  _{h} (V)(p)$ and $\tau  _{v} (V)(p)$ denote the vectors tangent to $M$ at $p$, such that $\tau  (V) _p = \{ \tau  _{h} (V)(p) \}^h  + \{ \tau  _{v} (V)(p)\} ^v $.
\end{Th}

\begin{Rk}
Since the condition $\tau (V)=0$ has a tensorial character, as usual we can assume it as a {\em definition} of  harmonic maps even when $M$ is not compact, and Theorem \ref{harmtens} extends at once to the non-compact case. 
\end{Rk}

\begin{Rk}
We now specify (\ref{htension}) and (\ref{vtension}) for classical metrics on $TM$. 

\medskip 
a) When $G=g^s$ is the Sasaki metric, we find the well known result:
$V:(M,g) \rightarrow (TM,g^s)$ is a harmonic map if and only if 
\begin{eqnarray}\label{tenssasaki}
& & {\rm tr} [R(\nabla _{\cdot} V,V)\cdot ]=0  \quad {\rm and}  \\
& & \bar \Delta V =0.
\end{eqnarray}

\medskip 
b) When $G=g_{CG}$ is the Cheeger-Gromoll metric, then, using (\ref{CheGro}), (\ref{htension}) and (\ref{vtension}), we easily get that $V:(M,g) \rightarrow (TM,g_{CG})$ is a harmonic map if and only if 
\begin{eqnarray}\label{tensCheGro}
& & {\rm tr} [R(\nabla _{\cdot} V,V)\cdot ]=0  \quad {\rm and}  \\
& & (1+r^2)\bar \Delta V + \nabla _{{\rm grad} \, r ^2} V
-\frac{1}{1+r^2}\left[ \vphantom{\displaystyle\frac{A}{A}} (2+r^2)
||\nabla V||^2 + \frac{1}{4}\left| \left| {\rm grad} \, r ^2
\right| \right|^2 \right] V =0.
\end{eqnarray}
Note that horizontal harmonicity of a vector field $V$, with respect to $g^s$ and $g_{CG}$, are expressed by the same condition.
\end{Rk}

We can now apply Theorem \ref{harmtens} to investigate relationships between harmonicity of maps defined by some special vector fields and properties of $g$-natural metrics.

\medskip\noindent
{\bf a) Parallel vector fields}

\noindent
It is well known that the existence of a non-vanishing parallel vector field $V$ on a Riemannian manifold $(M,g)$ is equivalent to the local reducibility of $M$ as $\Real \times M'$, equipped with the product metric, and $V$ is (locally) identified with a vector field tangent to the flat component $\Real$ of the product. Rewriting (\ref{htension}) and (\ref{vtension}) for a parallel vector field $V$, we have that $\tau (V)=0$ (and so, $V$ defines a harmonic map from $(M,g)$ to $(TM, G)$) if and only if 
\begin{equation}\label{htenspar}
-2 A_1 (\rho) QV + \left[ 2A_2 -A_3 g(QV,V)+nA_4  +\rho A_5 \right](\rho)  V =0 
\end{equation}
and 
\begin{equation}\label{vtenspar}
  -\bar \Delta V - B_1(\rho)  QV +\left[ 2B_3 -B_4 g(QV,V)+nB_5 + \rho B_6 \right] (\rho) V=0, 
 \end{equation}
where $\sqrt{\rho} = ||V||$ is the constant length of $V$.

Since $V$ is tangent to the flat component $\Real$ of the local decomposition $M=\Real \times M'$, it annihilates the curvature. Moreover, $\bar \Delta V =0$ for a parallel vector field. Therefore, (\ref{htenspar}) and (\ref{vtenspar})
are equivalent to 
\begin{equation}\label{tenspar}
\left[ 2A_2 +nA_4  +\rho A_5 \right](\rho)  = \left[ 2B_3 +nB_5 + \rho B_6 \right] (\rho) =0. 
\end{equation}
Using (\ref{Ai}) and (\ref{Bi}), we can easily conclude that (\ref{tenspar}) gives exactly (\ref{condpar}), that is,
$$[ (n-1)(\alpha _1 +\alpha _3 )+ \phi _1 +\phi _3 ]' (\rho)=0. $$
Therefore, we have the following

\begin{Th}\label{harmpar}
A parallel vector field $V$ defines a harmonic map $V:(M,g) \rightarrow (TM,G)$ if and only if its constant length satisfies {\rm (\ref{condpar})}, that is, $\rho = ||V|| ^2$ is a critical point of the function 
$$(n-1)(\alpha _1 +\alpha _3 )+ \phi _1 +\phi _3.$$
In particular:

\medskip
(i) For any Riemannian $g$-natural metric $G$ on $TM$ satisfying 
\begin{equation}\label{parharm}
 (n-1)(\alpha _1 +\alpha _3 )+ \phi _1 +\phi _3= {\rm constant},
\end{equation}
all parallel vector fields define harmonic maps from $(M,g)$ to $(TM, G)$.

\medskip
(ii) For any Riemannian $g$-natural metric $G$ on $TM$ such that 
\begin{equation}\label{noparharm}
[ (n-1)(\alpha _1 +\alpha _3 )+ \phi _1 +\phi _3 ]' (t) \neq 0 \; {\rm for \, all }\, t,
\end{equation}
parallel vector fields do not define harmonic maps from $(M,g)$ to $(TM, G)$.
\end{Th}

\noindent
Both $g^s$ and $g_{CG}$ satisfy (\ref{parharm}), as it easily follows from (\ref{Sasa}) and (\ref{CheGro}), respectively.  Hence, case {\em (i)} of Theorem \ref{harmpar} applies to both $g^s$ and $g_{CG}$.

\medskip\noindent
{\bf b) Vector fields of constant length}

\noindent
Consider a vector field $V \in \mathfrak X ^{\rho} (M)$. Then, from (\ref{htension}) and (\ref{vtension}) we get at once the following

\begin{Pro}\label{rhotens}
A vector field $V \in \mathfrak X ^{\rho} (M)$ satisfies $\tau (V)=0$ {\em (and so, it defines a harmonic map $V: (M,g) \rightarrow(TM,G)$)} if and only if
\begin{eqnarray}\label{unithtens}
 & & \quad -2 A_1 (\rho) QV +2 C_1 (\rho) {\rm tr} [R(\nabla _{\cdot} V,V)\cdot ] +2C_2 (\rho) \nabla _V V +\left[ 2A_2 (\rho)-A_3 (\rho) g(QV,V) \right.\\ 
 & &\left.  +nA_4 (\rho) +\rho A_5 (\rho) +2C_4 (\rho) g({\rm tr} [R(\nabla _{\cdot} V,V)\cdot ],V) +2C_5 (\rho) {\rm div}V +E_2 (\rho) ||\nabla V||^2 \right] V =0 \nonumber 
\end{eqnarray}
and 
\begin{eqnarray}\label{unitvtens}
& & -\bar \Delta V - B_1 (\rho) QV +2 D_1 (\rho) {\rm tr} [R(\nabla _{\cdot} V,V)\cdot ]
+2D_2 (\rho) \nabla _V V +\left[ 2B_3 (\rho) -B_4 (\rho) g(QV,V) \right.\\
 & & \left. +nB_5 (\rho)+ \rho B_6 (\rho)+2D_4 (\rho) g({\rm tr} [R(\nabla _{\cdot} V,V)\cdot ],V) +2D_5 (\rho) {\rm div}V   +F_2 (\rho) ||\nabla V||^2 \right]  V=0. \nonumber
 \end{eqnarray}
\end{Pro}

\noindent
In the special case when $(M,g)$ has constant sectional curvature $k$, from Proposition \ref{rhotens} it follows

\begin{Cor}\label{rhotenscostcurv}
Let $(M,g)$ be a Riemannian manifold of constant sectional curvature $k$. A vector field $V \in \mathfrak X ^{\rho} (M)$ defines a harmonic map from $V: (M,g) \rightarrow(TM,G)$) if and only if
\begin{eqnarray}\label{unithtensccurv}
& & 2 ( k C_1  +C_2)(\rho) \nabla _V V +\left[-2(n-1)k A_1 +2A_2 - (n-1)k \rho  A_3  \right.
\\ & & \left.  +nA_4  +\rho A_5  +2(C_5  -k C_1 - k \rho C_4 ) {\rm div}V +||\nabla V||^2 E_2 \right](\rho) V =0 \nonumber 
\end{eqnarray}
and 
\begin{eqnarray}\label{unitvtensccurv}
  \quad & & -\bar \Delta V + 2 ( k D_1  +D_2)(\rho) \nabla _V V +\left[-(n-1)k B_1 +2B_3 - (n-1)k \rho  B_4  \right.
\\ & & \left.  +nB_5  +\rho B_6  +2(D_5  -k D_1 - k \rho D_4 ) {\rm div}V +||\nabla V||^2 F_2 \right](\rho) V =0 \nonumber 
 \end{eqnarray}
\end{Cor}

\medskip\noindent 
{\bf Proof.} Since $(M,g)$ has constant sectional curvature $k$, its curvature tensor $R$ is given by
\begin{equation}\label{ccurv}
R(X,Y)Z=k(g(Y,Z)X- g(X,Z)Y).
\end{equation}
By (\ref{ccurv}) it easily follows that $QV=(n-1)k V$ and ${\rm tr} [R(\nabla _{\cdot} V,V)\cdot ]=-k ({\rm div}V )V+k \nabla _V V$. Using these formulae in (\ref{unithtens}) and (\ref{unitvtens}), we respectively get (\ref{unithtensccurv}) and (\ref{unitvtensccurv}) $\Box$

\bigskip\noindent
Equations (\ref{unithtens}) and (\ref{unitvtens}) are rather difficult to manage in full generality. We consider now the special case of a Riemannian $g$-natural metric $G$ for which $\alpha _2 (\rho)=\beta _2 (\rho)=0$. Note that $\alpha _2=\beta _2=0$ has a clear geometric meaning, since it characterizes $g$-natural metrics on $TM$ with respect to which horizontal and vertical distributions are mutually orthogonal. Under the assumption $\alpha _2 (\rho)=\beta _2 (\rho)=0$, taking into account formulae (\ref{Ai})-(\ref{Fi}) of the Levi-Civita connection of a $g$-natural metric $G$ given in Proposition \ref{lev-civ-con}, (\ref{unithtens}) and (\ref{unitvtens}) reduce respectively to
\begin{equation}\label{unitspech}
 C_1 (\rho) {\rm tr} [R(\nabla _{\cdot} V,V)\cdot ] +C_2 (\rho) \nabla _V V +\left[ C_4 (\rho) g({\rm tr} [R(\nabla _{\cdot} V,V)\cdot ],V) +C_5 (\rho) {\rm div}V \right] V =0 
\end{equation}
and 
\begin{equation}\label{unitspecv}
  -\bar \Delta V +\left( 2B_3 (\rho) +nB_5 (\rho)+\rho B_6 (\rho)+F_2 (\rho) ||\nabla V||^2 \right) V=0.
 \end{equation}
In particular, (\ref{unitspecv}) implies at once that $\bar \Delta V$ is collinear with $V$. So, $V$ is an eigenvector for the rough Laplacian $\bar \Delta$ and, since $\sqrt{\rho} = ||V||$ is a constant, we have $\bar \Delta V = \frac{1}{\rho} ||\nabla V||^2 V$ and (\ref{unitspecv}) implies
\begin{equation}\label{unitspecv1}
  \left( F_2 (\rho)-\frac{1}{\rho} \right) ||\nabla V||^2 +\left(2B_3 +nB_5 +t B_6 \right) (\rho) =0.
 \end{equation}
Again taking into account $\alpha _2 (\rho)=\beta _2 (\rho)=0$, (\ref{Bi}) and (\ref{Fi}), (\ref{unitspecv1}) may be easily rewritten as follows:
\begin{equation}\label{unitspecv2}
\left( \frac{1}{\rho} \alpha _1 + \alpha ' _1 \right) (\rho) ||\nabla V||^2 + \left[ (n-1) (\alpha _1 +\alpha _3) +\phi _1 +\phi _3 \right]' (\rho) =0.
 \end{equation}
Because of (\ref{unitspecv2}), for different Riemannian $g$-natural metrics $G$ some very different situations can occur about the harmonicity of the map $V:(M,g) \rightarrow (TM, G)$ defined by $V \in \mathfrak X ^{\rho} (M)$. The results are resumed in the following

\begin{Th}\label{unith}
Let $(M,g)$ be a Riemannian manifold and  $G$ a Riemannian $g$-natural metric on $TM$ satisfying $\alpha _2 (\rho)=\beta _2 (\rho)=0$,  $\rho >0$. Then, a vector field $V \in \mathfrak X ^{\rho} (M)$ defines a harmonic map $V:(M,g) \rightarrow (TM,G)$ if and only if it satisfies {\rm(\ref{unitspech})} and {\rm(\ref{unitspecv})}. In particular:

\medskip
(i) If  
\begin{equation}\label{unitspecv4}
\left( \frac{1}{\rho} \alpha _1 + \alpha ' _1 \right) (\rho) = \left[ (n-1) (\alpha _1 +\alpha _3) +\phi _1 +\phi _3 \right]' (\rho) =0,
 \end{equation}
then $V \in \mathfrak X ^{\rho} (M)$  defines a harmonic map $V:(M,g) \rightarrow (TM, \bar G)$ if and only if $V$ is an eigenvector of $\bar \Delta$ and {\rm(\ref{unitspech})} holds.

\medskip
(ii) If  
\begin{equation}\label{unitspec5}
\left( \frac{1}{\rho} \alpha _1 + \alpha ' _1 \right) (\rho) \neq 0 = \left[ (n-1) (\alpha _1 +\alpha _3) +\phi _1 +\phi _3 \right]' (\rho) ,
 \end{equation}
then $V \in \mathfrak X ^{\rho} (M)$ defines a harmonic map $V:(M,g)\rightarrow(TM, \bar G)$ if and only if $V$ is parallel. 

\medskip
(iii) If  
\begin{equation}\label{unitspec6}
\left( \frac{1}{\rho} \alpha _1 + \alpha ' _1 \right) (\rho) = 0 \neq \left[ (n-1) (\alpha _1 +\alpha _3) +\phi _1 +\phi _3 \right]' (\rho) ,
 \end{equation}
there are not vector fields $V \in \mathfrak X ^{\rho} (M)$ defining harmonic maps from  $(M,g)$ to $(TM, \bar G)$. 

\medskip
(iv) If  
\begin{equation}\label{unitspec7}
\left( \frac{1}{\rho} \alpha _1 + \alpha ' _1 \right) (\rho) \neq 0 \neq \left[ (n-1) (\alpha _1 +\alpha _3) +\phi _1 +\phi _3 \right]' (\rho) ,
 \end{equation}
then $V \in \mathfrak X ^{\rho} (M)$ defines a harmonic map $V: (M,g) \rightarrow (TM, G)$ if and only if {\rm(\ref{unitspech})} holds, $\bar \Delta V$ is collinear to $V$ and the length of $\nabla V$ satisfies
\begin{equation}\label{unitspec8}
\displaystyle  ||\nabla V|| ^2 = - \frac{ \rho \left[ (n-1) (\alpha _1 +\alpha _3) +\phi _1 +\phi _3 \right]' (\rho)}{ \left( \alpha _1 + \rho \, \alpha ' _1 \right) (\rho)}.
 \end{equation}
\end{Th}

\noindent
In the case of the Cheeger-Gromoll metric $g_{CG}$, (\ref{CheGro}) easily implies $\left( \frac{1}{\rho} \alpha _1 + \alpha ' _1 \right) (\rho) \neq 0 $ and $\left[ (n-1) (\alpha _1 +\alpha _3) +\phi _1 +\phi _3 \right]' (\rho) =0$. Therefore, {\em when $TM$ is equipped with $g_{CG}$, the parallel ones are the only vector fields of constant length, defining harmonic maps}. In particular, when $(M,g)$ has constant sectional curvature $k \neq 0$, then a vector field of constant length never defines a harmonic map $V:(M,g) \rightarrow (TM, g_{CG})$. 

\medskip
It is worthwhile to emphasize that, since a general Riemannian $g$-natural metric $G$ depends on six different smooth functions $\alpha _1$, $\alpha _2$, $\alpha _3$, $\beta _1$, $\beta _2$ and $\beta _3$ (satisfying inequalities (\ref{cond-riem})), in each of cases (i)-(iv) listed in Theorem \ref{unith}, there are plenty of Riemannian $g$-natural  metrics which furnish examples. We now illustrate some interesting cases:

\medskip\noindent
{\bf Example A:} {\em Assume $(M,g)$ has constant sectional curvature $k$. For any $\varepsilon >0$, there exists a family of Riemannian
$g$-natural metrics $\{ G_{\varepsilon} \}$, such that for all $\rho \geq \varepsilon$, $V \in \mathfrak X ^{\rho} (M)$ defines a harmonic map from $(M,g)$ to $(TM,G_{\varepsilon})$ if and only
if $\bar \Delta V$ is collinear to $V$.}

\medskip\noindent
In fact, it suffices to consider the family of $g$-natural metrics $\{ G_{\varepsilon} \}$ defined by the functions
\begin{equation}\label{ex1}
 \left\{
 \begin{array}{l}
   \alpha_1(t)=\lambda/t,  \quad (\textup{for} \; t \geq \varepsilon, \;
   \textup{and prolonged smoothly and positively to} [0,\varepsilon)),\\
\alpha_2=\beta_2=0, \\
  \alpha_1 +\alpha_3 =\mu , \\
\beta_1+\beta_3= -k\alpha_1, \\
\beta_1 \quad \textup{arbitrary such that} \; \alpha_1(t)
+t\beta_1(t)>0 \; \textup{for all} \; t>0,
  \end{array}
  \right.
\end{equation}
where $\lambda>0$ and $\mu>\sup(0,k\lambda)$. Formulae (\ref{ex1}) ensure that each $G_{\varepsilon}$ is Riemannian and, for all $\rho \geq \varepsilon$, we are in case (i) of Theorem \ref{unith}. Moreover, (\ref{unitspech}), equivalently (\ref{unithtensccurv}),  is satisfied. Note that whenever $\varepsilon \leq 1$, this case applies to Hopf vector fields of an odd-dimensional sphere.

\medskip\noindent
{\bf Example B:} {\em For any $\delta >0$, there exists a family of Riemannian $g$-natural metrics $\{ G_{\delta} \}$, such that for all $\rho \geq \delta$, $V \in \mathfrak X ^{\rho} (M)$ never defines a harmonic map from $(M,g)$ to $(TM,G_{\delta})$.}

\medskip\noindent
To show this, we consider the family of $g$-natural metrics $\{ G_{\delta} \}$ described by  
\begin{equation}\label{ex2}
  \left\{
  \begin{array}{l}
    \alpha_1(t)=\lambda/t, \quad (\textup{for} \; t \geq \varepsilon, \; \textup{and prolonged
smoothly to} \; [0,\epsilon)),\\
\alpha_2=\beta_2=0, \\
   \alpha_1 +\alpha_3 =\mu , \\
(\beta_1+\beta_3)(t)= \eta/t^2  \quad (\textup{for} \; t \geq \varepsilon, \; \textup{and
prolonged smoothly to} \; [0,\varepsilon)),\\
\beta_1 \quad \textup{arbitrary such that} \; \alpha_1(t)
+t\beta_1(t)>0 \; \textup{for all} \; t>0,
  \end{array}
  \right.
\end{equation}
for some positive constants $\lambda, \eta$. Then, each $G_{\delta}$ is Riemannian and for all $\rho \geq \delta$, we are in case (iii) of Theorem \ref{unith} $\Box$

\bigskip
As concerns the meaning of condition {\rm(\ref{unitspech})}, notice that, since $V \in \mathfrak X ^{\rho} (M)$, {\rm(\ref{unitspech})} implies
\begin{equation}\label{spec2}
(C_1 + \rho C_4)(\rho) g({\rm tr} [R(\nabla _{\cdot} V,V)\cdot ],V) + \rho C_5 (\rho) {\rm div}V =0. 
\end{equation}
When $M$ is compact, then $\int _M {\rm div} V dv_g =0$ and (\ref{spec2}) reduces to  
$$( C_1 + \rho C_4)(\rho)\int _M g({\rm tr} [R(\nabla _{\cdot} V,V)\cdot ],V) dv_g =0,$$
which in particular is satisfied whenever 
\begin{equation}\label{unitspec9}
{\rm tr} [R(\nabla _{\cdot} V,V)\cdot ]=0.
\end{equation}
Moreover, using formulae of Proposition \ref{lev-civ-con}, we can conclude that if $(\beta _1 +\beta _3)(\rho)=0$ (and $\alpha _2 (\rho)=\beta _2 (\rho)=0$), then
\begin{equation}\label{unispec9'}
 \left\{
\begin{array}{l}
C_1 (\rho) = -\frac{\alpha _1 ^2}{2 \alpha} >0, \\
C_2 (\rho) = C_4 (\rho) = C_5 (\rho) =0 \\
\end{array}
\right.
\end{equation}
and so, {\rm(\ref{unitspech})} reduces to (\ref{unitspec9}). We now apply this information to the special case of  Killing vector fields of constant length. 

The early theory of harmonic unit vector fields developped by Gil-Medrano and other authors (see \cite{G2} for a survey) shows that there are many interesting contexts in which  non-parallel unit vector fields satisfying (\ref{unitspec9}) appear. 

Let $V \in \mathfrak X(M)$ be a {\em Killing} vector field. As it is well-known, $V$ satisfies
\begin{equation}\label{spec1}
QV = \bar \Delta V.
\end{equation}
In the special case of an {\em Einstein} manifold $M$, we have $QV =  \frac{S}{n} V$, $S$ being the {\em scalar curvature} of $(M,g)$. Therefore, if $V \in \mathfrak X ^{\rho} (M)$, by (\ref{spec1}) it then follows
\begin{equation}\label{spec4}
\bar \Delta V = \displaystyle ||\nabla V||^2 V = \frac{S}{n} V.
\end{equation} 
 
Consider now any Riemannian $g$-natural metric $G$ on $TM$, satisfying $\alpha _2 (\rho)=\beta _2 (\rho)=0$ and 
\begin{equation}\label{spec5}
(t \alpha _1)' (\rho) \frac{S}{n} = - \rho \left[ (n-1) (\alpha _1 +\alpha _3) +\phi _1 +\phi _3 \right]' (\rho) .
 \end{equation}
Because of (\ref{spec4}) and (\ref{spec5}), we can conclude that (\ref{unitspecv}) and (\ref{unitspecv1}) (equivalently, (\ref{unitspecv2})) are satisfied. Therefore, if $\alpha _2 (\rho)=\beta _2 (\rho)=0$, a Killing vector $V \in \mathfrak X ^{\rho} (M)$ defines a harmonic map $V: (M,g) \rightarrow (TM,G)$ if and only if 
\begin{equation}\label{newspec}
C_1(\rho) {\rm tr} [R(\nabla _{\cdot} V,V)\cdot ] + C_4 (\rho) g({\rm tr} [R(\nabla _{\cdot} V,V)\cdot ],V)V =0. 
\end{equation}
Assuming also $(\beta _1 + \beta _3)(\rho) =0$, (\ref{unispec9'}) holds and hence, (\ref{newspec}) is equivalent to requiring that ${\rm tr} [R(\nabla _{\cdot} V,V)\cdot ]=0$. So, we have at once the following  

\begin{Th}\label{Killing}
Let $(M,g)$ be an Einstein manifold, $V \in \mathfrak X^{\rho} (M)$ a Killing vector field and $G$ a Riemannian $g$-natural metric on $TM$, satisfying 
$$
\left\{
\begin{array}{l}
\alpha _2 (\rho)=\beta _2 (\rho)=0, \\ (\beta _1 + \beta _3)(\rho) =0, \\
(t \alpha _1)' (\rho) \frac{S}{n} = - \rho \left[ (n-1) (\alpha _1 +\alpha _3) +\phi _1 +\phi _3 \right]' (\rho). 
\end{array}
\right.
$$
Then, $V$ defines a harmonic map $V:(M,g) \rightarrow (TM, G)$ if and only if ${\rm tr} [R(\nabla _{\cdot} V,V)\cdot ]=0$.
\end{Th}

\section {\red Riemannian $g$-natural metrics having parallel vector fields as the only harmonic sections} \label{ph}

Theorem \ref{harmpar} shows under which assumptions on a Riemannian $g$-natural metric $G$, a parallel vector field $V$ defines a harmonic map $V:(M,g) \rightarrow (TM,G)$. Obviously, this also includes the case of the Sasaki metric $g^s$ on $TM$. In fact, when $M$ is compact, parallel vector fields are all and the ones defining harmonic maps from $(M,g)$ into $(TM,g^s)$ \cite{I}, \cite{N}. We now consider the question whether this rigidity property is peculiar to the Sasaki metric, or there are other Riemannian $g$-natural metrics having the same property. 

By Theorem \ref{harmtens}, a vector field $V \in \mathfrak X (M)$ defines a harmonic map $V:(M,g) \rightarrow (TM, G)$ if and only if both (\ref{htension})
and (\ref{vtension}) are satisfied. Moreover, (\ref{parharm}) gives a necessary condition on these metrics, for the harmonicity of parallel vector fields. However, (\ref{parharm}) is not sufficient in general to conclude that a vector field $V$, satisfying (\ref{htension}) and (\ref{vtension}), is parallel. 

Looking for some special forms of equations (\ref{htension})
and (\ref{vtension}), we determine two classes of Riemannian $g$-natural metric $G$, for which harmonic sections are all and the ones parallel vector fields. The first class, which also includes the Sasaki metric as special case, was determined starting from the hypothesis $\alpha _2 =\beta _2=0$. The second class was found assuming that coefficients appearing in corresponding terms of (\ref{htension}) and (\ref{vtension}), are proportional, that is, there exists some constant $k \in \Real$ such that $2A_1 = k B_1$, $C_1 = k D_1$,....  This permits to remarkably simplify (\ref{htension}) and (\ref{vtension}). The results we found, with the complete description of the sets of conditions determining these two classes of $g$-natural metrics (conditions which also take into account of (\ref{cond-riem})), are resumed in the following

\begin{Th}\label{hp1}
Let $(M,g)$ be a compact Riemannian manifold and $G$ be a Riemannian $g$-natural metric on $TM$, satisfying one of the following sets of conditions:

\medskip\noindent
either 
\begin{equation}\label{fam1}
\left\{
\begin{array}{l}
\alpha _2 =\beta _2=0, \\
\alpha _1 = {\rm constant} >0, \\
\alpha _3 = {\rm constant} > -\alpha _1, \\
\beta _1 =-\beta _3  \geq 0, \\
\beta ' _1 \leq 0,
\end{array}
\right.
\end{equation}
or
\begin{equation}\label{fam2}
\left\{
\begin{array}{l}
\alpha _1 = {\rm constant} >0, \\
\alpha _2 = {\rm constant} \neq 0, \\
\alpha = \alpha _1 (\alpha _1 + \alpha _3)-\alpha ^2 _2 > 0, \\
\beta _1 =\beta _2 = 0,  \\
\beta _3 >0, \\
(n-1) (\alpha _1 +\alpha _3) +\phi _1 +\phi _3 = {\rm constant}.
\end{array}
\right.
\end{equation}

\medskip\noindent
Then, for any $V \in \mathfrak X(M)$, $V$ defines a harmonic map $V:(M,g) \rightarrow (TM,G)$ if and only if $V$ is parallel.
\end{Th}

\medskip\noindent
{\bf Proof.} We first notice that (\ref{fam1}) (or (\ref{fam2})) implies (\ref{cond-riem}). So, $g$-natural metrics described by (\ref{fam1}) (or (\ref{fam2})) are Riemannian.

Suppose now that (\ref{fam1}) holds. If $V \in \mathfrak X(M)$ is parallel, denote by $\rho$ the constant value of $||V||^2$. Because of Theorem \ref{harmpar}, harmonicity of $V: (M,g) \rightarrow (TM,G)$ is equivalent to (\ref{condpar}).
By (\ref{fam1}), $\alpha _1 +\alpha _3$ is constant and $\beta _1 + \beta _3 =0$. So, 
$$[ (n-1)(\alpha _1 +\alpha _3 )+ \phi _1 +\phi _3 ]'=[ n(\alpha _1 +\alpha _3 )+ \beta _1 +\beta _3 ]' =0. $$
Hence, (\ref{condpar}) is satisfied, that is, $V$ defines a harmonic map into $(TM,G)$.   

As concerns the converse, if $V$ defines a harmonic map into $(TM, G)$, then by Theorem \ref{harmtens}, (\ref{htension}) and (\ref{vtension}) hold. Starting from formulae of Proposition \ref{lev-civ-con} and using (\ref{fam1}), it is easy to check that  (\ref{vtension}) becomes
\begin{equation}\label{met2}
-\bar \Delta V + \left[ \frac{\beta _1(r^2) }{\phi _1(r^2)} ||\nabla V ||^2 + \frac{\beta ' _1(r^2)}{4 \phi _1(r^2)} || {\rm grad} \,r^2 || ^2 \right] V =0,
\end{equation}
where $r^2=||V||^2$. We take the scalar product of (\ref{met2}) by $V$ and integrate over $M$. Since $\int _M g(\bar \Delta V,V)dv_g=\int _M ||\nabla V||^2 dv_g$, taking into account the definition of $\phi _1$, we get
\begin{equation}\label{met3}
\int _M \frac{\alpha _1 (r^2)}{\phi _1 (r^2)} ||\nabla V||^2 dv_g - \int _M \frac{\beta ' _1 (r^2)}{4 \phi _1 (r^2)}|| {\rm grad} \,r^2 || dv_g =0. 
\end{equation}
By (\ref{cond-riem}) it follows $\alpha _1, \phi _1 >0$. Moreover, by (\ref{fam1}), $\beta _1 ' \leq 0$. Therefore, (\ref{met3})  
implies that $V$ is parallel.

Next, assume (\ref{fam2}) holds. Then (\ref{condpar}) is satisfied and so, a parallel vector field is harmonic. Conversely, let $V \in \mathfrak X (M)$ define a harmonic map into $(TM,G)$. By Theorem \ref{harmtens}, (\ref{htension}) and (\ref{vtension}) hold. Moreover, (\ref{fam2}) implies that there exists a constant $k \neq 0$ (more explicitly, $k=-\alpha _1/\alpha _2$), such that
\begin{eqnarray*}
& & (2A_1,C_1,C_3,E_1,C_2,A_2,A_3,A_4,A_5,C_4,C_5,C_6,E_2,E_3) \\ &=& k \, (B_1,D_1,D_3,F_1,D_2,B_3,B_4,B_5,B_6,D_4,D_5,D_6,F_2,F_3).
\end{eqnarray*}
We then divide (\ref{vtension}) by $k$ and substract (\ref{vtension}) by (\ref{htension}), and we obtain $\bar \Delta V=0$. Hence, $0=\int _M g(\bar \Delta V,V)dv_g=\int _M ||\nabla V||^2 dv_g$ and so, $V$ is parallel $\Box$

\bigskip
Note that (\ref{fam1}) determines a family of Riemannian $g$-natural metrics, depending on two real parameters $\alpha _1$ and $\alpha _3$ and a smooth function $\beta _1 : \Real ^+ \rightarrow \Real$ (satisfying some inequalities). Inside this class, the Sasaki metric  $g^s$ is the special case determined by $\alpha _1 =1$ and $\alpha _3 =\beta _1=0$. 

On the other hand, (\ref{fam2}) also determines a family of Riemannian $g$-natural metrics, depending on two real parameters $\alpha _1$ and $\alpha _2$, and a smooth function $\alpha _3: \Real ^+ \rightarrow \Real$, satisfying some inequalities. In fact, using the definitions of $\phi _i$, the last equation of (\ref{fam2}) permits to write down $\beta _3$ in function of $\alpha _1$, $\alpha _2$ and $\alpha _3$. Obviously, this class does not contain  the Sasaki metric, since in (\ref{fam2}) we must have $\alpha _2 \neq 0$.

\section {\red Critical points for the energy restricted to vector fields} \label{harmonicvf}

Let $(M,g)$ be a compact Riemannian manifold. We want to investigate conditions under which the map $V:(M,g) \rightarrow (TM,G)$ associated to a vector field $V \in \mathfrak X(M)$, is a critical point for the energy functional $E :\mathfrak X (M) \rightarrow \Real$, that is, only considering variations among maps defined by vector fields. Gil-Medrano \cite{G1} proved that, equipping $TM$ with the Sasaki metric $g^s$, $V:(M,g) \rightarrow (TM,g^s)$ is a critical point for the energy functional $E: \mathfrak X (M) \rightarrow \Real$ if and only if $V$ is parallel.

Consider now a vector field $V \in \mathfrak X(M)$, and a smooth variation $\{ V_t \} \subset \mathfrak X (M)$ of $V$, with $|t|< \varepsilon$ and $V_0 =V$. Note that $\pi  \circ V_t = {\rm id}_M$ for all $t$, where $\pi:TM \rightarrow M$ is the natural projection and ${\rm id}_M$ the identity on $M$. Therefore, the {\em variational vector field} $\widetilde W$ associated to the variation satisfies
$$ \widetilde W_p = \frac{\partial V_t (p)}{\partial t} \left|\begin{array}{l} \\ \hspace{-1mm} t=0 \end{array} \hspace{-2mm} \in \mathcal V _{V_p} TM, \right. $$
for all $p \in M$ and so, $V$ is a critical point for $E : \mathfrak X (M) \rightarrow \Real$ if and only if 
\begin{equation}\label{harm1}
0= E' (0) = \frac{d E_t}{d t} \left|\begin{array}{l} \\ \hspace{-1mm} t=0 \end{array} \hspace{-1mm} = - \int _M  G _{V_p} \left(\tau(V) _p,  \widetilde W _p \right) dv_g, \right. 
\end{equation}
for all variation $\{ V_t \} \subset \mathfrak X (M)$ of $V$. Note that, as was already remarked in \cite{G1}, for any vertical vector field $W^v$, section of the bundle $V^{-1} TTM$ of vector fields along $V$, there exists a variation $\{ V_t \} \subset \mathfrak X (M)$ of $V$, such that $\displaystyle W^v =\frac{\partial V_t }{\partial t} \left|\begin{array}{l} \\ \hspace{-1mm} t=0\end{array}\right.\hspace{-1mm}$. So, by (\ref{harm1}) it follows that $V$ is a critical point for $E : \mathfrak X (M) \rightarrow \Real$ if and only if 
\begin{equation}\label{harm2}
\int _M G _{V_p} \left(\tau(V) _p,  W^v _p \right) dv_g =0,
\end{equation}
for all vector fields $W \in \mathfrak X (M)$. Taking into account Proposition \ref{$g$-nat}, we easily find that (\ref{harm2}) is equivalent to 
\begin{eqnarray}\label{harm3}
& & \int _M g  \left( \vphantom{\displaystyle\frac{a}{a}} \alpha _2 \tau  _{h} (V) +  \beta _2 g (\tau  _{h} (V),  V ) V  +\alpha _1 \tau  _{v} (V) + \beta _1  g (\tau  _{v} (V),  V) V,  W \right)dv_g =0, 
\end{eqnarray}
where $\tau  _{h} (V)(p)$, $\tau  _{v} (V)(p)$ denote the vectors tangent to $M$ at $p$, such that $\tau  (V) _p = \{ \tau  _{h}(V)(p) \} ^h + \{ \tau  _{v} (V)(p)\} ^v$, for all $p \in M$. Since (\ref{harm3}) must hold for any vector field $W \in \mathfrak X (M)$, it is equivalent to requiring that
\begin{eqnarray}\label{harm4}
& & T(V) := \alpha _2 \tau  _{h}(V) +  \beta _2 g (\tau  _{h} (V),  V) V  +\alpha _1 \tau  _{v}(V) + \beta _1  g (\tau  _{v} (V),  V ) V=0. 
\end{eqnarray}
Note that (\ref{harm4}) expresses the {\em vanishing of the projection of the tension field $\tau (V)$ into the vertical distribution}, with respect to an arbitrary Riemannian $g$-natural metric $G$. 

Clearly, if $V:(M,g) \rightarrow (TM,G)$ is a harmonic map, in particular $V$ is a critical point for $E: \mathfrak X (M) \rightarrow \Real$. This is also expressed by formula (\ref{harm4}). In fact, if $V:(M,g) \rightarrow (TM,G)$ is a harmonic map, then Theorem \ref{harmtens} implies that $\tau _h (V)=\tau _v (V)=0$ and so, (\ref{harm4}) holds. In general, the converse does not hold. To emphasize this, we consider the special situation when $\alpha _2 = \beta _2 =0$. Under this assumption, $T(V)=0$ is equivalent to requiring that $\tau _v (V)=0$. In fact, taking into account $\alpha _2 = \beta _2 =0$, if $\tau _v (V)=0$ we have at once $T(V)=0$. Conversely, if $T(V)=0$, from $\alpha _2 = \beta _2 =0$ it follows that (\ref{harm4}) reduces to
\begin{eqnarray}\label{harm4spec}
& & \alpha _1 \tau  _{v}(V) + \beta _1  g (\tau  _{v} (V),  V ) V=0. 
\end{eqnarray}
Taking the scalar product of both sides of (\ref{harm4spec}) by $V$,
since $\phi _1 = \alpha _1 + r^2 \beta _1$, where $r ^2 = ||V||^2$,
we get
\begin{eqnarray}\label{harm4spec1}
& & \phi _1 g( \tau  _{v}(V), V ) =0.
\end{eqnarray}
By (\ref{cond-riem}), $\phi _1 >0$. Therefore, (\ref{harm4spec1}) gives $g( \tau  _{v}(V), V ) =0$ and (\ref{harm4spec}) reduces to $\alpha _1 \tau  _{v}(V)=0$. Again (\ref{cond-riem}) gives $\alpha _1 >0$ and so, $\tau _v (V)=0$. In this way, by Theorem \ref{harmtens} we obtain at once the following 

\begin{Th}\label{char}
Let $(M,g)$ be a compact Riemannian manifold and $G$ any
Riemannian $g$-natural metric on $TM$ with 
$\alpha_2=\beta_2=0$. A vector field $V$ on $M$ defines a harmonic
map $V:(M,g) \rightarrow (TM,G)$ if and only if the following
conditions hold:
\begin{itemize}
    \item [i)] $\tau_h(V)=0$, and
    \item [ii)] $V$ is a critical point for $E: \mathfrak X (M) \rightarrow \Real$, that is, $T(V)=0$.
\end{itemize}    
\end{Th}

Coming back to the general case, we recall that formulae (\ref{htension}), (\ref{vtension}) describe the tension field associated to $V:(M,g) \rightarrow (TM, G)$, for an arbitrary  Riemannian $g$-natural metric $G$. Using (\ref{htension}), (\ref{vtension}) in (\ref{harm4}) and taking into account (\ref{Ai})-(\ref{Fi}), some long but standard calculations lead to the following characterization:

\begin{Th}\label{harmvf}
Let $(M,g)$ be a compact Riemannian manifold. A vector field $V \in \mathfrak X (M)$ is a critical point for $E :\mathfrak X (M) \rightarrow \Real$ if and only if $T(V)=0$, where
\begin{eqnarray}\label{vfharm}  
\qquad  T(V) &=& -\alpha _1 \bar \Delta V  + \left( \alpha ' _2 - \frac{\beta _2}{2} \right) {\rm grad} \, r^2 +\alpha ' _1  \nabla _{{\rm grad} \, r ^2} V \\ 
 & &-\left\{ \vphantom{\displaystyle\frac{A}{A}} [(n-1)(\alpha _1 +\alpha _3) + \phi _1 + \phi _3]' + (2\alpha ' _2 -\beta _2) {\rm div}V + (\alpha ' _1 -\beta _1) ||\nabla V||^2    \right. \nonumber \\
& &\vphantom{\displaystyle\frac{a}{A}} \beta _1 g(\bar \Delta V,V) - \left( \phi _2 C_6 + \phi _1 D_6 + \beta _2 C_2 + \beta _1 D_2 + \beta _2 C_3 + \beta _1  D_3 \right)  V(r^2) \nonumber \\
& & \left. - \frac{1}{4} \left( \phi _2 E_3 + \phi _1 F_3 + 2\beta _2 E_1 + 2\beta _1 F_1 \right) \left| \left| {\rm grad} \, r ^2 \right| \right|^2   \right\} V, \nonumber 
\end{eqnarray}
and all fuctions are evaluated at $r^2 = ||V||^2$. 
\end{Th}

Since the critical point condition $T(V)=0$ has a tensorial character, it also makes sense when $(M,g)$ is not compact. For a general Riemannian manifold $(M,g)$, if a vector field $V$ satisfies $T(V)=0$, we call it a {\em $\mathfrak X$-harmonic} vector field. 


\begin{Rk}
Specifying (\ref{vfharm}) for the Sasaki and Cheeger-Gromoll metrics of $TM$, we have the following results:

\begin{itemize}
\item if $G=g^s$, then (\ref{Sasa}) implies the well known formula  $T(V)= - \bar \Delta V$. 
\item if $G=g_{CG}$, then applying (\ref{CheGro}) we easily obtain
\begin{eqnarray*}
T(V) &=& -\frac{1}{1+r^2} \bar \Delta V  -\frac{1}{(1+r^2)^2} \nabla _{{\rm grad} \, r ^2} V  \\
& & +\frac{1}{1+r^2} \left[ -g(\bar \Delta V,V) +  \frac{2+r^2}{1+r^2} ||\nabla V||^2 - \frac{1}{4(1+r^2)} \left| \left| {\rm grad} \, r ^2 \right| \right|^2 \right] V. \nonumber
\end{eqnarray*}
\end{itemize}
\end{Rk}

We now determine $\mathfrak X$-harmonic vector fields, under some special assumptions either on the vector fields themselves or on the Riemannian $g$-natural metric $G$.

\medskip
{\bf 1) Parallel vector fields.}

\noindent
Suppose $V \in \mathfrak X (M)$ is a parallel vector field. Then, $\nabla V=0$, $\bar \Delta V=0$ and $||V||^2=\rho$ is a constant. Thus, (\ref{vfharm}) reduces to
$$
T(V) = -[(n-1)(\alpha _1 +\alpha _3) + \phi _1 + \phi _3]' (\rho) \, V.
$$
Hence, $T(V)=0$ coincides with the necessary and sufficient condition we found in Theorem \ref{harmpar} for the harmonicity of $V : (M,g) \rightarrow (TM, G)$, and in Section 3 for critical points of the energy $E$ restricted to parallel vector fields. Therefore, we get the following 

\begin{Th}\label{par}
Let $(M,g)$ be a Riemannian manifold and $G$ any Riemannian $g$-natural metric on $TM$. For a parallel vector field $V$ on $M$, the following statements are equivalent:

\medskip
(a) $V : (M,g) \rightarrow (TM, G)$ is a harmonic map;

\medskip
(b) $V$ is $\mathfrak X$-harmonic;

\medskip
(c) $V$ is a critical point for $E$ in the set $\mathfrak X_{\mathcal P} (M)$ of all parallel vector fields on $M$;

\medskip
(d) $\rho =||V||^2$  is a critical point for the function $[(n-1)(\alpha _1 +\alpha _3) + \phi _1 + \phi _3]$.%

\end{Th}	
	
\noindent	
Theorem \ref{par} includes as special cases both the Sasaki metric $g^s$ and the Cheeger-Gromoll metric $g_{CG}$ on $TM$, for which (d) is trivially satisfied and so, all parallel vector fields define harmonic maps.

\medskip\goodbreak
{\bf 2) Vector fields of constant length.}

\noindent
Considering a vector field $V \in \mathfrak X ^{\rho} (M)$, by (\ref{vfharm}) we have that $T(V)=0$ if and only if 
\begin{eqnarray}\label{harm7}
& &\; \alpha _1 (\rho) \bar \Delta  V + \left\{ \beta _1 (\rho) g(\bar \Delta V,V)+ [(n-1)(\alpha _1 +\alpha _3) + \phi _1 + \phi _3]'(\rho)  \right.\\
 & & \left. + (2 \alpha ' _2 - \beta _2 )(\rho)  {\rm div}V + (\alpha ' _1 -\beta _1)(\rho) ||\nabla V||^2  \right\}V=0. \nonumber
\end{eqnarray}
By (\ref{harm7}) it follows at once that $\bar \Delta V$ is collinear to $V$. Therefore, since $V$ has constant length $||V||=\sqrt{\rho}$, we have $\bar \Delta V = \frac{1}{\rho}||\nabla V||^2 V$ and from (\ref{harm7}) we get 
\begin{eqnarray}\label{harm9}
& & \left( \frac{1}{\rho} \alpha _1 + \alpha ' _1 \right)(\rho) ||\nabla V||^2 +(2 \alpha ' _2 - \beta _2 ) (\rho) {\rm div}V  + [(n-1)(\alpha _1 +\alpha _3) + \phi _1 + \phi _3]' (\rho)  =0. 
\end{eqnarray}
Thus, Theorem \ref{harmvf} implies the following

\begin{Th}\label{uni}
Let $(M,g)$ be a Riemannian manifold and $G$ any Riemannian $g$-natural metric on $TM$. A vector field $V \in \mathfrak X ^{\rho} (M)$ is  $\mathfrak X$-harmonic if and only  if $\bar \Delta V$ is collinear to $V$ and {\em (\ref{harm9})} holds.
\end{Th}	
 
\noindent
For the Sasaki metric $g^s$, an arbitrary vector field $V$ is $\mathfrak X$-harmonic if and only if $\nabla V=0$. For the Cheeger-Gromoll metric $g_{CG}$, by Theorem \ref{uni} we have the following 

\begin{Cor}\label{CheGrocrit}
Let $(M,g)$ be a Riemannian manifold and equip $TM$ with the Cheeger-Gromoll metric $g_{CG}$. A vector field $V \in \mathfrak X ^{\rho} (M)$ is $\mathfrak X$-harmonic if and only  if it is parallel {\em (and so, if and only if $V:(M,g)\rightarrow(TM,g_{CG})$ is harmonic)}.
\end{Cor}	

\medskip\noindent 
{\bf Proof.} If $V \in \mathfrak X (M)$ is parallel, then the conclusion follows from Theorem \ref{par}. Conversely, assume $V \in \mathfrak X ^{\rho} (M)$. Using (\ref{CheGro}), (\ref{harm9}) gives at once $\nabla V=0$ $\Box$

\bigskip
Equation (\ref{harm9}) remains quite difficult to solve in full generality. For this reason, we consider the special case when $\alpha _2 (\rho)=\beta _2 (\rho)=0$. Under this assumption, (\ref{harm9}) becomes
\begin{eqnarray}\label{harm10}
& & 2 \alpha ' _2 (\rho) {\rm div}V + \left( \frac{1}{\rho} \alpha _1 + \alpha ' _1 \right)(\rho) ||\nabla V||^2 + [(n-1)(\alpha _1 +\alpha _3) + \phi _1 + \phi _3]' (\rho)  =0. 
\end{eqnarray}
In particular, if $\alpha ' _2 (\rho)=0$, then (\ref{harm10}) gives exactly (\ref{unitspecv2}) which, together with the collinearity of $\bar \Delta V$ and $V$, is equivalent to (\ref{unitspecv}). Therefore, calculations above, together with Theorem \ref{unith}, lead at once to the following

\begin{Pro}\label{unispec}
Let $(M,g)$ be a Riemannian manifold and $V \in \mathfrak X ^{\rho} (M)$. For any  Riemannian $g$-natural metric $G$ on $TM$, satisfying $\alpha _2 (\rho)=\alpha ' _2 (\rho)=\beta _2 (\rho)=0$,

\medskip
(1) $V$ is $\mathfrak X$-harmonic if and only if {\em (\ref{unitspecv})} holds.

(2) $V$ defines a harmonic map $V:(M,g) \rightarrow (TM, G)$ if and only if it is $\mathfrak X$-harmonic and satisfies {\em (\ref{unitspech})}. 

\medskip\noindent
In particular, $\mathfrak X$-harmonic vector fields do not necessarily define harmonic maps. 
\end{Pro}	

\noindent	
Taking into account formulae (\ref{ex1}) determining the Riemannian $g$-natural metrics given in Example A, from Proposition \ref{unispec} we obtain the following

\begin{Cor}
Let $(M,g)$ be a Riemannian manifold of constant sectional curvature $k$. For any $\varepsilon >0$, there exists a family of Riemannian $g$-natural
metrics $\{ G_{\varepsilon} \}$, such that for all $\rho \geq
\varepsilon$, $V \in \mathfrak X ^{\rho} (M)$ defines a harmonic
map from $(M,g)$ to $(TM,G_{\varepsilon})$ if and only if it is $\mathfrak X$-harmonic.
\end{Cor}

\section {\red Harmonicity of the Reeb vector field} \label{expl}

We now apply the previous study to the case of some classic vector fields, namely, Reeb vector fields and Hopf vector fields, and we start by recalling some basic definitions and properties about contact metric manifolds. 

Given a smooth manifold $M$ of odd dimension $n=2m+1$, a {\em contact structure} $(\eta , \varphi
,\xi)$ over $M$ is composed by a global $1$-form $\eta$ (the {\em contact form}) such that
$\eta \wedge (d \eta) ^{m} \neq 0$ everywhere on $M$, a
global vector field $\xi$ (the {\em Reeb {\em or} characteristic vector field})
and a global tensor $\varphi $, of type (1,1), such that
\begin{equation}\label{reeb0}
    \eta (\xi) =1 \,, \quad \varphi \xi =0 \,, \quad \eta \varphi =0
    \,, \quad \varphi ^2 = -I+ \eta \otimes \xi \,.
\end{equation}
A Riemannian metric $g$ is said to be {\em associated} to the
contact structure $(\eta , \varphi ,\xi)$, if it satisfies
\begin{equation}\label{reeb1}
    \eta = g(\xi , \cdot)\,, \quad d \eta =g(\cdot ,\varphi \cdot) \,,
    \quad g(\cdot, \varphi \cdot) = -g(\varphi \cdot ,\cdot ) \,.
\end{equation}    
We refer to $(M, \eta ,g)$ or to $(M, \eta ,g,\xi ,
\varphi)$ as a contact metric manifold. As it is well known, the Reeb vector field $\xi$ plays a very important role in describing the geometry of a contact metric manifold. By (\ref{reeb0}) and (\ref{reeb1}) it follows at once that $\xi$ is a unit vector field on $(M,g)$, that is, $\xi \in \mathfrak X ^1 (M)$.

As it is well-known, the Reeb vector field $\xi$ satisfies
\begin{equation}\label{reeb2}
     \nabla \xi = -\varphi - \varphi h ,
     \quad \nabla _{\xi} \xi =0, \quad {\rm div} \xi = 0, 
\end{equation}
where $h = \frac 12  \cal L _\xi \varphi$ is the Lie derivative of $\varphi$, and 
\begin{equation}\label{reeb3}
    || \nabla \xi ||^2 =2m+ {\rm tr} h^2 = 4m -g(Q \xi,\xi).
\end{equation}
Moreover, as it was proved  in \cite{P2}, 
\begin{equation}\label{reeb4}
 \bar \Delta \xi =4m \xi -Q \xi .
\end{equation}
For further details, references and information about contact metric manifolds, we refer to \cite{B}.

In \cite{P2}, the third author introduced and studied {\em $H$-contact spaces}, that is, contact metric manifolds $(M, \eta ,g,\xi ,\varphi)$  whose Reeb vector field $\xi$ is a critical point for the energy functional $E$ restricted to the space $\mathfrak X ^{1}(M)$ of all unit vector fields on $(M,g)$, considered as smooth maps from $(M,g)$ into the {\em unit tangent sphere bundle} $T^1 M$, equipped with the Riemannian metric induced on $T^1 M$ by the Sasaki metric $g^s$ of $TM$. As it was proved in \cite{P2}, $(M, \eta ,g,\xi ,\varphi)$ is $H$-contact if and only if $\xi$ is an eigenvector of the Ricci operator. (As a consequence, the class of $H$-contact manifolds is very large, since $\eta$-Einstein spaces, $K$-contact spaces, $(k,\mu)$-spaces and strongly locally $\phi$-symmetric spaces are all $H$-contact.) 

We now use (\ref{reeb2})-(\ref{reeb4}) to rewrite (\ref{unithtens}) and (\ref{unitvtens}) for $\xi$. By Proposition \ref{rhotens}, we then get the following
\begin{Pro}\label{xiharm}
Let $(M, \eta ,g,\xi ,\varphi)$ be a contact metric manifold and $G$ an arbitrary Riemannian $g$-natural metric on $TM$. The Reeb vector field $\xi$ defines a harmonic map $\xi: (M,g) \rightarrow(TM,G)$ if and only if
\begin{eqnarray}\label{xihtens}
& & -2 A_1 (1) Q \xi +2 C_1 (1) {\rm tr} [R (\nabla _{\cdot} \xi, \xi)\cdot ]+\left[ \vphantom{\displaystyle\frac{a}{a}} 2A_2 (1)+(2m+1)A_4 (1)+A_5 (1)  \right.
\\ & & \left. + 4m E_2 (1) +2C_4 (1) g({\rm tr} [R (\nabla _{\cdot} \xi, \xi)\cdot ],\xi) -[A_3 (1) +E_2 (1)] g(Q\xi,\xi) \vphantom{\displaystyle\frac{a}{a}} \right] \xi =0 \nonumber 
\end{eqnarray}
and 
\begin{eqnarray}\label{xivtens}
  & &[1- B_1 (1)] Q\xi +2 D_1 (1) {\rm tr} [R (\nabla _{\cdot} \xi, \xi)\cdot ]
+\left[ \vphantom{\displaystyle\frac{a}{a}} -4m + 2B_3 (1) +(2m+1)B_5 (1)  \right.
\\ & & \left. +B_6 (1) +4m F_2 (1)+2D_4 (1) g({\rm tr} [R (\nabla _{\cdot} \xi, \xi)\cdot ],\xi) -[B_4 (1)+F_2 (1)] g(Q\xi,\xi) \vphantom{\displaystyle\frac{a}{a}} \right]  \xi=0. \nonumber
 \end{eqnarray}
\end{Pro}
Since $C_1=-\frac{\alpha_1 ^2}{2\alpha} \neq 0$, we can use (\ref{xihtens}) to write $ {\rm tr} [R (\nabla _{\cdot} \xi, \xi)\cdot ]$ as a linear combination of $Q\xi$ and $\xi$, and we get 
\begin{eqnarray}\label{trR}
& &{\rm tr} [R (\nabla _{\cdot} \xi, \xi)\cdot ] =\frac{1}{2C_1(1)} \left\{  2 A_1 (1) Q \xi  -\left[ \vphantom{\displaystyle\frac{a}{a}} 2A_2 (1)+(2m+1)A_4 (1)+A_5 (1)  \right. \right.
\\ & &\left. \left. + 4m E_2 (1) +2C_4 (1) g({\rm tr} [R (\nabla _{\cdot} \xi, \xi)\cdot ] -[A_3 (1) +E_2 (1)] g(Q\xi,\xi) \vphantom{\displaystyle\frac{a}{a}} \right] \xi \right\}. \nonumber 
\end{eqnarray}
Replacing (\ref{trR}) in (\ref{xivtens}), we obtain
\begin{eqnarray}\label{Hcont}
  & &\left[ 1- B_1 (1)+\frac{2A_1(1)D_1(1)}{C_1(1)}\right] Q\xi + \left\{ \left[ \vphantom{\displaystyle\frac{a}{a}} -4m
  + 2B_3 (1) +(2m+1)B_5 (1) +B_6 (1) \right. \right.
\\ & & \left. +4m F_2 (1) +2D_4 (1) g({\rm tr} [R (\nabla _{\cdot} \xi, \xi)\cdot ],\xi)-[B_4 (1)+F_2 (1)] g(Q\xi,\xi) \vphantom{\displaystyle\frac{a}{a}} \right]  \nonumber \\
& &-\frac{D_1(1)}{C_1(1)} \left[ \vphantom{\displaystyle\frac{a}{a}} 2A_2 (1)+(2m+1)A_4 (1)+A_5 (1) + 4m E_2 (1) +2C_4 (1) g({\rm tr} [R (\nabla _{\cdot} \xi, \xi)\cdot ],\xi)  \right.\nonumber
\\ & &\left. \left. -[A_3 (1) +E_2 (1)] g(Q\xi,\xi) \vphantom{\displaystyle\frac{a}{a}} \right] \right\} \xi  =0. \nonumber 
\end{eqnarray}
Note that, by (\ref{Ai})-(\ref{Di}), we easily see that
$$B_1 =\frac{2A_1D_1}{C_1} .$$
So, (\ref{Hcont}) implies that $\xi$ is a Ricci eigenvector and we have at once the following

\begin{Th}\label{xiarm}
Let $(M, \eta ,g,\xi ,\varphi)$ be a contact metric manifold and $G$ an arbitrary Riemannian $g$-natural metric on $TM$. If $\xi$ defines a harmonic map $\xi: (M,g) \rightarrow(TM,G)$, then $(M,\eta,g)$ is $H$-contact.
\end{Th}

Under some assumptions on the Riemannian $g$-natural metric $G$, we are able to completely characterize harmonicity of $\xi : (M,g) \rightarrow (TM,G)$. In particular, if $\alpha _2 (1) =\beta _2 (1)=0$, then (\ref{xihtens}) and (\ref{xivtens})
reduce to
\begin{equation}\label{xih}
 C_1 (1) \, {\rm tr} [R (\nabla _{\cdot} \xi, \xi)\cdot ] + C_4 (1)\, g({\rm tr} [R (\nabla _{\cdot} \xi, \xi)\cdot ],\xi) \xi =0 
\end{equation}
and 
\begin{equation}\label{xiv}
  Q \xi= \left[ 4m- 2B_3 (1) -(2m+1)B_5 (1)-B_6 (1)- 4m F_2 (1) +F_2 (1) g(Q \xi,\xi) \right] \xi,
 \end{equation}
respectively. (\ref{xiv}) means that $\xi$ is a Ricci eigenvector, that is, $M$ is $H$-contact. Moreover, by (\ref{xiv}), the corresponding Ricci eigenvalue $g(Q \xi,\xi)$ depends on functions which determine the metric $G$. On the other hand, by (\ref{reeb3}) we have $g(Q\xi,\xi) = 2m -{\rm tr} h^2$ and so, (\ref{xiv}) is equivalent to requiring that $Q \xi$ is collinear to $\xi$ and
\begin{equation}\label{xiv'}
 [F_2 (1)-1] {\rm tr} h^2 = - \left[ 2B_3 (1) +(2m+1)B_5 (1)+B_6 (1)+ 2m F_2 (1)-2m \right] .
 \end{equation}
Notice that since $\xi$ is a unit vector, (\ref{xiv'}) also follows from (\ref{unitspecv}).
Taking into account formulae (\ref{Ai})-(\ref{Fi}), we can write coefficients of (\ref{xiv'}) explicitly in function of $\alpha _i$, $\beta _i$. Thus, (\ref{xiv'}) becomes
\begin{equation}\label{newxiv}
 \,( {\rm tr} h^2 + 2m )\,(\alpha _1 +\alpha ' _1)(1)  + [2m(\alpha _1+\alpha _3) + \phi _1+\phi _3]'(1)=0 .
 \end{equation}

As concerns (\ref{xih}), note that taking the scalar product of (\ref{xih}) by $\xi$ and by an arbitrary vector field $X$ orthogonal to $\xi$, we obtain
\begin{equation}\label{xih'}
\left\{
\begin{array}{l}
[C_1 (1)+C_4 (1)]\,  g({\rm tr} [R (\nabla _{\cdot} \xi, \xi)\cdot ],\xi)=0, \vphantom{\displaystyle\frac{a}{a}} \\ 
 C_1 (1) \, g({\rm tr} [R (\nabla _{\cdot} \xi, \xi)\cdot ], X)=0 \; {\rm for \, all \,} X \perp \xi. \vphantom{\displaystyle\frac{a}{a}}
\end{array}
\right.
\end{equation}
As we already noticed, $C_1=-\frac{\alpha_1 ^2}{2\alpha} \neq 0$. Moreover, since $\alpha _2 (1) =\beta _2 (1)=0$, by
(\ref{cond-riem}) and the definition of $\phi$, $\phi _i$ we easily get
$$ C_1 (1) +C_4 (1)=-\frac{\alpha_1}{2(\phi_1+\phi_3)} \neq 0.$$ 
Because of (\ref{xih'}), (\ref{xih}) is equivalent to requiring  ${\rm tr} [R( \nabla _{\cdot} \xi, \xi)\cdot ]=0$.
So, we have the following

\begin{Th}\label{xicrispec}
Let $(M, \eta ,g,\xi ,\varphi)$ be a contact metric manifold and $G$ any Riemannian $g$-natural metric on $TM$, satisfying $\alpha _2 (1) =\beta _2 (1)=0$. Then $\xi$ defines a harmonic map $\xi: (M,g) \rightarrow(TM,G)$ if and only if  $M$ is $H$-contact, {\em (\ref{newxiv})} holds and ${\rm tr} [R( \nabla _{\cdot} \xi, \xi)\cdot ]=0$.
\end{Th}

\begin{Rk}
a)\, We recall that a unit vector field $U$ defines a harmonic map $U: (M,g) \rightarrow(T^1M,g^s)$ if and only if 
$\bar\Delta U$ is collinear to $U$ and ${\rm tr} [R( \nabla _{\cdot} U, U)\cdot ]=0$ (see \cite{HYi}). For a Riemannian $g$-natural metric $G$ on $TM$, satisfying $\alpha _2 (1) =\beta _2 (1)=0$ and (\ref{newxiv}), Theorem \ref{xicrispec} gives the following interesting fact: {\em $\xi :(M,g) \rightarrow (TM,G)$  is  a harmonic map if and only if $\xi : (M,g) \rightarrow (T^1M,g^s)$ is a harmonic map}.\\
b)\, When $(\alpha_1+\alpha ' _1)(1)=0$, formula (\ref{newxiv}) reduces to (\ref{condpar}). On the other hand, if $(\alpha_1+\alpha ' _1)(1)\neq 0$, then (\ref{newxiv})  implies that tr$h^2$ is constant. All homogeneous contact metric manifolds provide examples of contact metric spaces for which tr$h^2$ is a constant. 
\end{Rk}

A {\em $K$-contact space} is a contact metric manifold $(M, \eta ,g,\xi ,\varphi)$ satisfying $h=0$. As it was remarked in \cite{P2}, a $K$-contact space is necessarily $H$-contact. For a $K$-contact space, (\ref{newxiv}) clearly reduces to  
\begin{equation}\label{Kcont}
2m (\alpha _1 +\alpha ' _1)(1) +[2m(\alpha _1+\alpha _3) + \phi _1+\phi _3]'(1)=0.
\end{equation}

\noindent
As it is well-known, {\em Sasakian manifolds} are $K$-contact, while the converse only holds in dimension three. Assume now $(M, \eta ,g,\xi ,\varphi)$ is Sasakian and consider a {\em $\varphi$-basis} on $M$, that is, an orthonormal basis of vector fields  $\{ e_1,..,e_m,\varphi e_1,..,\varphi e_m, \xi \}$. Taking into account (\ref{reeb0}), the first equation in (\ref{reeb2}) and the first Bianchi identity, we can see that the Reeb vector field $\xi$ satisfies
\begin{eqnarray*}
-{\rm tr} [R( \nabla _{\cdot} \xi, \xi)\cdot ] &=& {\rm tr} [R( (\varphi + \varphi h) \cdot, \xi)\cdot ]=\sum _{i=1} ^m 
[R(\varphi e_i,\xi)e_i+R(\varphi ^2 e_i,\xi) \varphi e_i] \\ 
&=& \sum _{i=1} ^m [R(\varphi e_i,\xi)e_i -R(e_i,\xi) \varphi e_i] = -\sum _{i=1} ^m 
R(e_i, \varphi e_i)\xi =0, 
\end{eqnarray*}
since on a Sasakian manifold, $R(X,Y)\xi=0$ for all $X,Y$ orthogonal to $\xi$ \cite{B}. Hence, Theorem \ref{xicrispec} implies the following

\begin{Th}\label{Sas}
Let $(M, \eta ,g,\xi ,\varphi)$ be a Sasakian manifold, dim$M=2m+1$ and $G$ any Riemannian $g$-natural metric on $TM$, satisfying $\alpha _2 (1) =\beta _2 (1)=0$. Then, $\xi$ defines a harmonic map $\xi: (M,g) \rightarrow(TM,G)$ if and only if {\em (\ref{Kcont})} holds.
\end{Th}

Next, we shall investigate under which conditions the Reeb vector field is $\mathfrak X$-harmonic. Since $\xi$ is a unit vector field, it is $\mathfrak X$-harmonic  if and only if (\ref{harm7}) holds. Moreover, taking into account (\ref{reeb3}) and (\ref{reeb4}), (\ref{harm7}) becomes
\begin{eqnarray}\label{xicrit}
 & & \alpha _1 (1) Q \xi = \left\{ 4m(\alpha _1 +\alpha ' _1)(1) + [2m(\alpha _1 +\alpha _3) + (\phi _1 + \phi _3)]'(1)   - \alpha ' _1 (1) g(Q\xi,\xi)  \right\}\xi. 
\end{eqnarray}
Since $\alpha _1 >0$, (\ref{xicrit}) gives that $\xi$ is a Ricci eigenvector. Using this fact and (\ref{reeb3}), (\ref{xicrit}) reduces to (\ref{newxiv}).
Hence, from Theorem \ref{xicrispec} we obtain the following

\begin{Th}\label{xicri}
Let $(M, \eta ,g,\xi ,\varphi)$ be a contact metric manifold and $G$ an arbitrary Riemannian $g$-natural metric on $TM$. If $\xi$ is  $\mathfrak X$-harmonic, then $M$ is $H$-contact. Conversely, if $M$ is $H$-contact, then $\xi$ is  $\mathfrak X$-harmonic if and only if {\em (\ref{newxiv})} holds.
\end{Th}

\begin{Rk}
Note that (\ref{newxiv}) is not fulfilled neither by the Sasaki metric nor by the Cheeger-Gromoll metric on $TM$, as it is easy follows from (\ref{Sasa}) and (\ref{CheGro}), respectively. So, when $(M, \eta ,g,\xi ,\varphi)$ is an arbitrary contact metric manifold and $TM$ is equipped with either $g^s$ or $g_{CG}$, then the Reeb vector field $\xi$ is never $\mathfrak X$-harmonic. In particular, in such cases, $\xi$ never defines a harmonic map.

On the other hand, it is easy to exhibit examples of Riemannian $g$-natural metrics, satisfying (\ref{newxiv}). For example, (\ref{newxiv}) holds for all Riemannian $g$-natural metrics belonging to the two-parameters family  satisfying 
$$
\left\{
\begin{array}{l}
\alpha _1 (t)= k_1 e^{-t}, \\
\alpha _3 (t)=  k_2 -\alpha _1 (t), \\
\alpha _2 = \beta _1 =\beta _2 = \beta _3 =0,
\end{array}
\right.
$$
where $k_1,k_2$ are positive constants. 

\end{Rk}

\medskip\noindent
We now apply Theorem \ref{xicri} to special classes of contact metric manifolds, namely, $K$-contact and $(k,\mu)$-spaces. If we assume $(M, \eta ,g,\xi ,\varphi)$ is $K$-contact, then $Q \xi = 2m \xi$ and (\ref{xicrit}) becomes (\ref{Kcont}).

Next, we recall that a contact metric manifold $(M, \eta ,g,\xi ,\varphi)$ for which $\xi$ belongs to the $(\kappa,\mu)$-nullity
distribution, that is,
\begin{eqnarray}\label{kmu}
& &  R(X,Y)\xi = \kappa\big(\eta(Y)X -\eta(X)Y\big) + \mu\big(\eta(Y)h X -\eta(X)hY\big) , 
\end{eqnarray}
where  $\kappa,\mu$ are constants, is called a  {\em $(\kappa,\mu)$-space}. Such class of spaces extends that of Sasakian manifolds. The constant $\kappa$ satisfies $\kappa \leq 1$; if $\kappa =1$, then $\mu=0$ and $M$ is Sasakian (\cite{B}, Theorem 7.7). Moreover,  $(\kappa,\mu)$-spaces are examples of  strongly pseudo-convex CR manifolds (\cite{B}, Theorem 7.6), and non-Sasakian $(\kappa,\mu)$-spaces are examples of locally $\phi$-symmetric spaces (\cite{B}, p. 118). From (\ref{kmu}), one gets $Q \xi = 2m \kappa\,\xi$ and so, (\ref{xicrit}) becomes 
 \begin{eqnarray}\label{crikmu}
 \, 2m (2-\kappa )\,(\alpha _1 +\alpha ' _1)(1)  + [2m(\alpha _1+\alpha _3) + \phi _1+\phi _3]'(1)=0. 
\end{eqnarray}
Then, by Theorem \ref{xicri} we have the following

\begin{Th}\label{xicrisp}
Let $(M, \eta ,g,\xi ,\varphi)$ be a contact metric  manifold and $G$ any Riemannian $g$-natural metric on $TM$. 

\medskip
\, \ (i)\, If $M$ is $K$-contact, then
$\xi$ is $\mathfrak X$-harmonic if and only if {\em (\ref{Kcont})} 
holds. 

\medskip
\, \ (j)\, If $M$ is a $(\kappa,\mu)$-space, then $\xi$ is $\mathfrak X$-harmonic if and only if { \em(\ref{crikmu})} 
holds.
\end{Th}

\medskip
We now recall that {\em Hopf vector fields} on the unit sphere $S^{2m+1}$, equipped with its canonical metric $g_0$, are all and the ones Killing unit vector fields on $S^{2m+1}$ \cite{W2}. Moreover, a Hopf vector field $\bar \xi $ can always be considered as the Reeb vector field of a suitable Sasakian structure $(S^{2m+1}, \bar \eta, g_o,\bar \xi, \bar \varphi)$, where $\bar \eta =g_0 (\cdot, \bar \xi)$ and $\bar \varphi=-\nabla \bar \xi$. Taking into account Theorems \ref{Sas} and \ref{xicrisp} above, we have

\begin{Cor}\label{Hopfharm}
For all Riemannian $g$-natural metrics on $T S^{2m+1}$, satisfying $\alpha _2 (1)=\beta _2 (1)=0$, a Hopf vector field $\bar \xi$ defines a harmonic map  $\bar \xi:(S^{2m+1},g_0) \rightarrow (T S^{2m+1},G)$ if and only if {\em (\ref{Kcont})} holds.
\end{Cor}

\begin{Cor}\label{Hopfcrit}
For all Riemannian $g$-natural metrics on $T S^{2m+1}$,
a Hopf vector field $\bar \xi$ is $\mathfrak X$-harmonic if and only if {\em (\ref{Kcont})} holds.
\end{Cor}

\begin{flushleft}
\textsc{D\'epartement des Math\'ematiques, Facult\'e des sciences Dhar El Mahraz,\\
Universit\'e Sidi Mohamed Ben Abdallah, B.P. 1796, F\`es-Atlas, F\`es, Morocco.}\\
\textit{E-mail address}: mtk{\_}abbassi@Yahoo.fr.
\end{flushleft}

\begin{flushleft}
\textsc{Dipartimento di Matematica "E. De Giorgi", Universit\`{a} degli Studi di Lecce, Lecce, ITALY.} \\
\textit{E-mail address}: giovanni.calvaruso@unile.it, domenico.perrone@unile.it.
\end{flushleft}

\end{document}